\documentclass[12pt]{amsart}
\usepackage{amscd,amssymb,graphicx}
\author[Fialowski]{Alice Fialowski}
\address{
Alice Fialowski\\
University of P\'ecs and
E\"otv\"os Lor\'and University\\ Hungary
}
\email{fialowsk@ttk.pte.hu, fialowsk@cs.elte.hu}
\author{Michael Penkava}
\address{University of Wisconsin\\
Eau Claire, WI 54702-4004} \email{penkavmr@uwec.edu}
\thanks{Research of the authors was partially supported by
grants from the University of Wisconsin-Eau Claire.}

\theoremstyle{definition}

\def \ph{\varphi}

\def \diag{\operatorname {diag}}

\def \ra{\rightarrow}

\def \hom{\mbox{\rm Hom}}

\def \gl{\mbox{$\mathfrak{gl}$}}
\def \tns{\otimes}

\def \k{\mbox{$\mathbb K$}}

\def \C{\mbox{$\mathbb C$}}
\def \Z{\mbox{$\mathbb Z$}}

\def\zt{\mbox{$\Z_2$}}

\def\A{\mbox{$\mathcal A$}}

\def\m{\mbox{$\mathfrak m$}}
\def\a{\mbox{$\mathfrak a$}}

\def\linf{\mbox{$L_\infty$}}
\def\and{\mbox{ \rm and }}

\def\s#1{(-1)^{#1}}

\def\pha#1#2{\ph^{#1}_{#2}}

\def\psa#1#2{\psi^{#1}_{#2}}

\def\P{\mathbb P}

\def\GL{{\mathbf{GL}}}

\def\sl{\mathfrak{sl}}

\begin{document}
\setlength{\multlinegap}{0pt}
\title{Moduli Spaces of low dimensional Lie superalgebras}%

\date{\today}
\begin{abstract}
In this paper, we study moduli spaces of low dimensional complex Lie superalgebras.
We discover a similar pattern for the structure of these moduli spaces as we observed
for ordinary Lie algebras, namely, that there is a stratification of the moduli space
by projective orbifolds. The moduli spaces consist of some families as well as some
singleton elements. The different strata are linked by jump deformations, which gives a
unique manner of decomposing the moduli space which is consistent with deformation theory.
\end{abstract}
\maketitle

\section{Introduction}
In a series of papers, the authors and some collaborators have been studying moduli spaces of low dimensional Lie algebras,
as well as moduli spaces of complex associative algebras, including algebras defined on a \zt-graded space.  In the Lie algebra case,
we have studied moduli spaces of complex Lie algebras of dimension up to 5, and real Lie algebras of dimension up to 4. In all of these
cases, we found that the moduli space has a natural decomposition into strata which are parameterized by projective orbifolds of a 
very simple kind, which
is a new point of view that had not appeared in the earlier literature.  

Each stratum was of the form $\P^n/G$, where $G$ is a subgroup of the symmetric group $\Sigma_{n+1}$, which acts on $\P^n$ by permuting the
projective coordinates. This led us to conjecture that every moduli space of finite dimensional Lie algebras has such a decomposition,
and it also reasonable to guess that the same conjecture holds for Lie superalgebras. In this paper we prove that the conjecture holds
for low dimensional complex Lie superalgebras

The classification of moduli
spaces of superalgebras is complicated by the fact that a Levi decomposition of a superalgebra does not always exist, and the fact
that for Lie superalgebras, a semisimple algebra may not be a direct sum of simple algebras.  However, the definitions of solvable
and nilpotent Lie superalgebras are the same as in the ordinary case.  A semisimple Lie algebra is one whose solvable radical (maximal
solvable ideal) is trivial.  Moreover, if a Lie superalgebra is not solvable, then the quotient by the solvable radical is semisimple.

If a Lie superalgebra $L$ is solvable, then it has a codimension 1 ideal, so there is an exact sequence
$$0\ra M\ra L\ra W\ra 0,$$
where $M$ is a \zt-graded ideal, and $W$ is a 1-dimensional algebra (which is necessarily trivial).  However, $W$ may be $1|0$ or $0|1$-dimensional.
As a consequence, solvable Lie superalgebras of a fixed dimension $m|n$ can be constructed from solvable Lie algebras of dimension $m-1|n$ or
dimension $m|n-1$, so this method allows one to construct the solvable $m|n$-dimensional Lie superalgebras by a bootstrap analysis.

The situation with superalgebras which are not solvable is more complex, but in low dimensions this complication mostly disappears owing to the
fact that there are not many examples of low dimensional complex semisimple Lie superalgebras.  The paper \cite{kac} contains a description
of the finite dimensional simple Lie superalgebras, as well as a prescription for constructing semisimple superalgebras.  A more recent article
\cite{hurni}, gives a more explicit description of the semisimple Lie superalgebras. The reader may also find the book \cite{scheunert} useful.

When a Lie superalgebra is not solvable, one has an exact sequence of the form

$$0\ra M\ra L\ra W\ra 0,$$ where this time, $W$ is semisimple and $M$ is the solvable radical.
Thus, for both solvable and nonsolvable Lie superalgebras, we can classify non semisimple algebras as extensions of either semisimple or trivial
Lie algebras by solvable Lie algebras.  There is a long history of the study of extensions of Lie algebras.  In \cite{fp9}, a description of the
process was given in the language of codifferentials, and the methods described in that article were used to construct the moduli spaces we are
studying here.

The bidimension of a \zt-graded vector space is given in the form $m|n$, where $m$ is the dimension of the even part of the space and $n$ is the dimension
of the odd part.  An ordinary $m$-dimensional Lie algebra is simply a Lie superalgebra of dimension $m|0$,  and it is necessary to consider such algebras
in studying the moduli spaces of superalgebras, because the dimension of the space $M$ or $W$ in the decomposition as an extension may be of the form $k|0$.
In fact, in the study of $3|1$-dimensional superalgebras, one has to consider the extension of the $3|0$-dimensional simple Lie algebra $\sl{(2,\C)}$ by
a $0|1$-dimensional (trivial) algebra.  However, other than this case, we won't have to consider any semisimple superalgebras in the moduli spaces
we construct,  owing to the low dimensions of the spaces involved.

In this paper, we address the complete moduli spaces for Lie superalgebras in dimensions $1|1$, $1|2$, $2|1$. $1|3$, $2|2$ and $3|1$.
\section{The language of codifferentials}
Classically, the space of cochains $C(L,L)$ of a Lie algebra with coefficients in the adjoint representation is given by $C(L,L)=\hom(\bigwedge L,L)$, where
$\bigwedge L$ is the exterior algebra of $L$.  For an ordinary algebra, $\bigwedge L$ has dimension $n2^n$, where $n=\dim(L)$, with components
$C^k(L,L)=\hom(\bigwedge^k L,L)$ of dimension $\binom nk$. There is a
natural Lie superalgebra structure on $C(L,L)$ and the Lie algebra structure itself is represented as an odd element $\ell$ of $C^2(L,L)$, which satisfies
the codifferential condition $[\ell,\ell]=0$.  It is possible to extend this definition to the \zt-graded case, but there is a more fundamental approach,
based on the fact that the exterior algebra of a \zt-graded space coincides in a natural manner with the symmetric algebra on the parity reversion of
the \zt-graded space.  Under this association, we obtain an equality between $C(L,L)$ and $C(\Pi L,\Pi L)=\hom(S(\Pi L),\Pi L)$. In fact, there is an isomorphism
between $C(\Pi L,\Pi L)$ and the space of \zt-graded coderivations of the symmetric coalgebra $S(\Pi L)$.  The difference in the expression of a Lie algebra
structure on $L$ and a codifferential on $S(L)$ is easy to express.  If $d$ is the codifferential on $S(W)$ corresponding to a Lie superalgebra structure $\ell$ on
$L$, then
$$\ell(a\wedge b)=\s{a}\pi(d(\pi a\cdot\pi b)).$$
Thus to convert from a Lie superalgebra expressed as a codifferential back to the standard form involves only  multiplication by a sign.
We will give our algebras in the form of codifferentials, but we will also indicate in some cases how to translate to the standard form.
\section{Construction of moduli spaces by extensions}
Let us assume that $0\ra M\ra L\ra W\ra 0$ gives an extension of the algebra structure on $W$ given by a codifferential $\delta$ by an
algebra structure on $M$ given by the codifferential $\mu$. Then if $d$ is the corresponding codifferential of the algebra structure on $L$,
we have $d=\delta+\mu+\lambda+\psi$, where $\lambda\in\hom(M\tns W,M)$ and $\psi\in\hom(S^2(W),M)$.  The term $\lambda$ is traditional called
the \emph{module} structure on $M$ and the term $\psi$ is called the \emph{cocycle}, although when $\mu\ne0$, $\lambda$ is not precisely
a module structure on $M$, nor is $\psi$ really a cocycle.  However, we will use this terminology (even though it is not precisely correct).
The condition that $d$ is a codifferential on $L$ is that $[d,d]=0$, which is equivalent to the following three conditions:
\begin{enumerate}
\item $[\mu,\lambda]=0$ (The compatibility condition)
\item $[\delta,\lambda]+\tfrac12[\lambda,\lambda]+[\mu,\psi]=0$ (The Maurer-Cartan Condition)
\item $[\delta+\lambda,\psi]=0$ (The cocycle condition)
\end{enumerate}
To construct an extension, we first fix $\delta$ and $\mu$, and then we solve the compatibility condition, which puts some constraints on
the coefficients of $\lambda$.  If $\beta\in\hom(W,M)$ is even, then replacing $\lambda$ with $\lambda+[\mu,\beta]$ generates an equivalent
extension, so we use this to simplify the form of $\lambda$.  Next, we consider the action of the group $G_{\delta,\mu}$ of transformations
of $M\oplus W$, consisting of those block diagonal elements such that the action of the appropriate piece on $\delta$ or $\mu$ preserves
this structure.  This allows us to restrict the form of $\lambda$ even more.

Next, we apply the Maurer-Cartan (MC) condition to $\lambda$ and a generic $\psi$, which may place additional constraints on the coefficients
of $\lambda$ and constraints on the coefficients of $\psi$.  Finally, we apply the cocycle condition, to construct a $d$ which is a codifferential.
Now, we can also apply a group $G_{\delta,\mu,\lambda}$ to restrict the coefficients of $\psi$ further, but in practice, we mostly did not
do this, except possibly for a diagonal transformation.

After doing this, we find some codifferentials, and study their equivalence classes.  Some of them naturally arise as families, and in this
case, we check to see that they represent a projective family, in the sense that multiplying the coefficients by a nonzero complex number
yields an isomorphic algebra.

This is not quite all the details involved, but we will illustrate the situation in our examples.

\section{Deformations of algebras and the versal deformation}

Given a 1-parameter family $d_t$ of algebras such that $d_0=d$, then we say that this family determines a deformation of $d$. If $d_t\not\sim d$ for $t$ in some
punctured nbd of $t=0$, then we say that the deformation is nontrivial.  If $d_t\sim d'$ for all nonzero $t$ in some punctured nbd of $0$, then this deformation
is called a \emph{jump deformation} of $d$, while if $d_s\not\sim d_t$ for $s\ne t$ for small enough $s$ and $t$, then  the deformation is called
a \emph{smooth deformation}.  One can also have multiparameter deformations $d_t$ where $t=(t_1,t_2,\cdots)$, in which case there may be 1-parameter curves
in the $t$ space which determine jump deformations and other curves which determine smooth deformations.

There is a generalization of these multiparameter deformations called a deformation with a \emph{local base} (see \cite{fi1}), which is a commutative algebra $A$ such that
there is an $A$-Lie algebra structure $d_A$ defined on $V\tns A$, where $V$ is the underlying vector space on which the Lie algebra is defined, and
$A$ is a local algebra, meaning it has a unique maximal ideal $\m$. One requires that $A/\m=\k$, the underlying field on which the Lie algebra is defined, so that there is a natural decomposition $A=\k\oplus\m$. Then there is a natural map $V\tns\A\ra V$, determined by the projection $A\ra\k$, and $d_A$ is called
a deformation with base $A$ if the induced map takes $d_A$ to $d$.

For super Lie algebras, it makes sense to work with \zt-graded commutative algebras, but
we don't take that point of view here, even though it would be interesting. If the reader
 is interested in seeing this type of analysis, we mention that in the study of low dimensional \linf\ algebras in \cite{fp5,fp6} we did consider
this more general perspective.

There is a special type of multiparameter deformation called a \emph{versal deformation}, has the property that it induces every deformation with a local
base in a natural manner.  Moreover, there is a special type of versal deformation, called a \emph{miniversal deformation} (see \cite{fi1}) which can be constructed in
a concrete fashion by beginning with an \emph{infinitesimal deformation} $d_1=d+t_i\delta^i$, where $\langle \overline{\delta}_1\rangle$ is a basis for $H^2(d)$, the second cohomology of the algebra $d$, see \cite{fifu}. The deformation is called infinitesimal because it satisfies the Jacobi identity up to first order terms in the $t_i$.
When studying Lie superalgebras, we only look at the even part of $H^2$, because we aren't considering deformations with a graded commutative base.

In \cite{fp1}, a method of constructing a miniversal deformationfor \linf\ algebras was outlined, and we have developed tools using the Maple computer algebra system for
carrying out the computations, which are mostly just applications of linear algebra, although  the computation of the versal deformation involves solving
systems of quadratic equations.

\section{The Moduli Space of $1|1$-dimensional Lie Superalgebras}
For ordinary $2$-dimensional Lie algebras, the moduli space consists of one nontrivial element, $\ell=\pha{e_1,e_2}{e_2}$, which is solvable, but not nilpotent.
Expressed as a codifferential, this solvable algebra has the form $d=\psa{1,2}2$.

For $1|1$-dimensional Lie superalgebras, the situation is more interesting. There are 2 nonequivalent $1|1$-dimensional Lie superalgebras.
Let $L=\rangle f,e\langle $ be a $1|1$-dimensional vector space with $e$ an even and $f$ an odd basis element.
The first algebra $\ell_1$ is given by $\ell_1(e,f)=f$. This algebra is analogous to the ordinary Lie algebra $\ell$ above. The second algebra
$\ell_2$, is given by the formula $\ell_2(f,f)=e$, with all other brackets vanishing.  Because $f$ is odd, $f\wedge f$ is not equal to zero, a
situation that cannot arise in the nongraded case.

The first algebra arises as an extension of the trivial algebra structure $\delta=0$ on a $0|1$-dimensional vector space $W=\langle v_2\rangle$ by
the trivial algebra structure $\mu=0$ on a $1|0$-dimensional space $M=\langle v_1\rangle$.  The structure $\lambda$ is determined up to
a constant multiple of $\psa{1,2}1$.  The structure $\psi$ must vanish as $S^2(W)=0$, since $W$ is a 1-dimensional odd vector space. Thus, up to isomorphism,
we obtain that the only possible nontrivial structure is $d_1=\psa{1,2}1$.

The second algebra arises as an extension of the trivial algebra structure $\delta=0$ on a $1|0$-dimensional vector space $W=\langle v_1\rangle$ by the
trivial algebra structure $\mu=0$ on a $0|1$-dimensional space $M=\langle v_2\rangle$. In the language of codifferential, all of the maps
$\delta$,$\mu$, $\lambda$ and $\psi$ must be odd, which forces $\lambda=0$, and $\psi$ to be a multiple of $\psa{1,1}2$.  Thus the only nontrivial structure
is given (up to isomorphism) by $d_2=\psa{1,1}2$.

One can also proceed directly to construct the algebras by noting that the form of the algebra must be $d=\psa{12}1a+\psa{1,1}2b$, and then checking that
the condition $[d,d]=0$ is equivalent to $ab=0$.

In the table below, we compute the cohomology of the two algebras. Here, $h_n$ is the bi-dimension  $H^n$, the cohomology of the algebra in degree $n$.
Let us recall the meaning of the cohomology in low degrees.  $H^0$ is the center of the algebra, $H^1$ is the space of nontrivial derivations of the algebra, $H^2$ classifies the infinitesimal deformations, and $H^3$ gives  the obstructions to extending an infinitesimal deformation.  An algebra for which the cohomology
vanishes in all degrees is called totally rigid.  In terms of actual deformations, only the odd part of $H^2$ counts, although one can interpret the even part in terms of deformations with a base given by a \zt-graded algebra.

The algebra $d_1$ is totally rigid.  We see that $H^0=\langle v_2\rangle$ is the center of the algebra $d_2$.  There is also a nontrivial even derivation of
$d_2$, given by $\pha{1}1+2\pha{2}2$.  This completely describes the moduli space of $1|1$-dimensional Lie superalgebras. Since $H^2=0$ for both of these algebras, there are no nontrivial deformations for either one of them.  Note that $d_1$ is solvable but not nilpotent, while $d_2$ is nilpotent.

\begin{table}[ht]
\begin{center}
\begin{tabular}{rlllrrrr}
&Algebra&&Codifferential&$h_0$&$h_1$&$h_2$&$h_3$\\
\hline\\
&$d_{1}$&=&$\psa{1,2}1$&$0|0$&$0|0$&$0|0$&$0|0$\\
&$d_{2}$&=&$\psa{1,1}2$&$0|1$&$1|0$&$0|0$&$0|0$\\
\hline
\end{tabular}
\end{center}
\caption{\protect Cohomology of $1|1$-Dimensional Complex Lie
Algebras}\label{Table 1}
\end{table}
\section{The moduli space of $2|1$-dimensional Lie Superalgebras}

There are only three (families of) algebras on a $2|1$-dimensional vector space.  The corresponding dimension for the codifferentials is $1|2$. One of these
is a projective family $d_3(p:q)$,  which means that $d_3(up:uq)\sim d_3(p:q)$ for all nonzero $u\in\C$.  In this case, there are no isomorphisms between
$d_3(p: q)$ and $d_3(x: y)$, except for the isomorphisms that give rise to the projective description in our notation.  As is usually the case,
when there is a family of algebras, there are some special values of the parameters $(p:q)$ such that the cohomology or even the deformation picture is
different than generically.  There is a special element $(0: 0)$, which is called somewhat unfortunately the generic element of projective space by
algebraic geometers, because the algebra corresponding to $(0:0)$ is never generic in its behavior.  In this case, $d_3(0:0)$ is actually the trivial
algebra, which has jump deformations to every nontrivial algebra in the moduli space. By a jump deformation, we mean a deformation $d_t$ of an algebra
$d$, where $t$ is a (multi)-index such that $d_t\sim d'$ for some algebra $d'$ except when $t=0$, in which case we obtain the original algebra.

The algebras $d_1$ and $d_3(p:q)$, except for $d_3(0: 0)$, are solvable but not nilpotent, while the algebras $d_2$ and $d_3(0:0)$ are nilpotent.
The algebra $d_2$ has a jump deformation to $d_1$, while $d_3(p: q)$ has a smooth deformation along the family. 
For the special cases of the parameters, only $d_3(1: -2)$ does not behave generically in terms of its deformations, because it has a jump deformation to
$d_1$ in addition to smooth deformations along the family.

\begin{table}[ht]
\begin{center}
\begin{tabular}{rlllrrrr}
&Algebra&&Codifferential&$h_0$&$h_1$&$h_2$&$h_3$\\
\hline\\
&$d_{1}$&=&$4\psa{1,1}2+\psa{1,3}1-2\psa{2,3}2$&$0|0$&$0|0$&$0|0$&$0|0$\\
&$d_{2}$&=&$4\psa{1,1}2$&$0|2$&$3|1$&$1|1$&$0|0$\\
&$d_{3}(p:q)$&=&$p\psa{1,3}1+q\psa{2,3}2$&$0|0$&$1|0$&$0|1$&$0|0$\\
&$d_{3}(1:-3)$&=&$\psa{1,3}1-3\psa{2,3}2$&$0|0$&$1|0$&$0|1$&$0|1$\\
&$d_{3}(1:0)$&=&$\psa{1,3}1$&$0|1$&$2|0$&$0|1$&$0|0$\\
&$d_{3}(0:1)$&=&$\psa{2,3}2$&$1|0$&$1|1$&$1|1$&$1|1$\\
&$d_{3}(1:-2)$&=&$\psa{1,3}1-2\psa{2,3}2$&$0|0$&$1|0$&$0|2$&$1|0$\\
&$d_{3}(0:0)$&=&$0$&$1|2$&$5|4$&$6|6$&$6|6$\\
\hline
\end{tabular}
\end{center}
\caption{\protect Cohomology of $2|1$-Dimensional Complex Lie
Algebras}\label{Table 2}
\end{table}

All of the $2|1$-dimensional algebras are given by extending the trivial algebra structure on either a $1|0$ dimensional algebra by an algebra structure on a $1|1$-dimensional space or by extending the trivial $0|1$-dimensional algebra by an algebra structure on a $2|0$-dimensional space. Because this case is fairly
easy to describe, but shows some of the important features of the construction, we will give an explicit construction of the moduli space.

First, let us consider $W=\langle v_3\rangle$, and $M=\langle v_1,v_2\rangle$, so that $v_1$ is the only even basis element. The module structure
$\lambda$ is of the form $\lambda=\psa{1,3}1a_1+\psa{2,3}2a_2$.  The $\psi$ term must vanish, and the map $\beta:W\ra M$ is of the form $\beta=\pha{3}2b$.

We know that $\mu$ is one of the $1|1$-dimensional algebras.  Let us first consider the case $\mu=\psi{1,2}1$. The compatibility condition forces $a_1=0$,
and then in this case, $\lambda=[\mu,\beta]$ for $b=a_2$, so we can assume that $\lambda$ vanishes.  As a consequence, all of the conditions for
$d$ to be an algebra structure are automatically satisfied, and we obtain only the algebra $d=\psa{1,2}1$, which is isomorphic to $d_3(1:0)$.

Next, let us assume that $\mu=\psa{22}1$. Then the compatibility condition gives $a_2=-2a_1$. Since $[\mu,\beta]=0$ for all $b$, we cannot simplify
the expression for $\lambda$, but applying the group $g_{\delta,\mu}$ we find we can reduce to the case $a_1=1$ or $a_1=0$.  The first case gives
the algebra $d_1$, while the second gives the algebra $d_2$.

Finally, we have to consider the case $\mu=0$. In this case, the compatibility condition is trivial, so we can express $\lambda=\psa{1,3}1p+\psa{2,3}2q$.
Here we substituted $p=a_1$ and $q=a_2$, which is a notation we use when we suspect the relation between the $p$ and the $q$ gives a projective symmetry.
In fact, the group $G_{\delta,\mu}$ acts on $\lambda$ by multiplying both coordinates by the same number, which is precisely what we expect if the
structures form a projective family.  We obtain $d=\lambda=d_3(p:q)$.

The reader may notice that we have already discovered all of the algebras by only looking at one of the possible decompositions. In fact, we suspect
that if $V$ is a solvable $n|1$-dimensional algebra with $n>1$, then there is a $(n-1)|1$-dimensional ideal. We already know from the study of
$1|1$-dimensional algebras  that the result does not hold for $n=1$.

However, let us proceed with the case $W=\langle v_2,v_3\rangle$ and $M=\langle v_1\rangle$.  First, we note that $\lambda$ and $\beta$ must both
vanish, while $\psi=\psa{1,1}2c_1+\psa{1,1}3c_2$.

Let us first consider the case when $\mu=\psa{1,2}1$.  In this case, we only need to consider the MC equation, and this forces $\psi=0$.  Thus
$d=\mu$ which is isomorphic to $d_3(0:1)$.

Finally, consider the case when $\mu=0$. Then there are no conditions on $\psi$, so we obtain $d=\psa{1,1}2c_1+\psa{1,1}3c_2$.  A little work with
$G_{\delta,\mu}$ would show that we only have to consider 2 cases, where $c_1=1$ and $c_2=0$, or both of the coefficients vanish.  In the first case,
$d\sim d_2$ while the second case is given by $d=d_3(0: 0)$, the trivial codifferential.

\section{The moduli space of $1|2$-dimensional Lie Superalgebras}
There is one projective family $d_1(p:q)$, and three singletons in this moduli space.
They correspond to codifferentials on a $2|1$-dimensional space.
Because the family $d_1(p:q)$ occurs first in the order of the description (we have ordered our algebras
in such a manner that an algebra only deforms to one whose number is lower, except for the generic element),
we don't have any room for extra deformations for special values of the parameter $(p:q)$ in $d_1(p:q)$.
There is an action of the symmetric group $\Sigma_2$ on the family by permuting the coordinates, which means
that $d_1(p:q)\sim d_1(q:p)$ for all $(p:q)$. This means that the family $d_1(p:q)$ is parametrized by
the projective orbifold $\P^1/\Sigma_2$.  This is a typical pattern that has arisen in our studies of moduli
spaces.

Note that the odd part of the dimension of $H^2$ is always 1, except for $d_1(0:0)$, where this number is 2.
It is always the case that the generic element in a family has jump deformations to every other element in the
family, and these are the only deformations of $d_1(0:0)$.

The element $d_2$ has a jump deformation to $d_1(1:1)$, while the other elements are rigid. The algebras
$d_1(0:0)$, $d_3$ and $d_4$ are nilpotent, while $d_2$ and $d_1(p:q)$ are solvable but not nilpotent otherwise.

\begin{table}[ht]
\begin{center}
\begin{tabular}{rlllrrrr}
&Algebra&&Codifferential&$h_0$&$h_1$&$h_2$&$h_3$\\
\hline\\
&$d_{1}(p:q)$&=&$p\psa{1,3}1+\psa{2,3}1+q\psa{2,3}2$&$0|0$&$1|0$&$0|1$&$0|0$\\
&$d_{1}(1:2)$&=&$\psa{1,3}1+\psa{2,3}1+2\psa{2,3}2$&$0|0$&$1|0$&$1|1$&$0|1$\\
&$d_{1}(1:3)$&=&$\psa{1,3}1+\psa{2,3}1+3\psa{2,3}2$&$0|0$&$1|0$&$0|1$&$1|0$\\
&$d_{1}(1:-1)$&=&$\psa{1,3}1+\psa{2,3}1-\psa{2,3}2$&$0|0$&$1|0$&$0|1$&$1|0$\\
&$d_{1}(1:0)$&=&$\psa{1,3}1+\psa{2,3}1$&$1|0$&$1|1$&$1|1$&$1|1$\\
&$d_{1}(0:0)$&=&$\psa{2,3}1$&$1|0$&$2|1$&$2|2$&$2|2$\\
&$d_{2}$&=&$\psa{1,3}1+\psa{2,3}2$&$0|0$&$3|0$&$0|3$&$0|0$\\
&$d_{3}$&=&$\psa{1,2}3+2\psa{1,1}3$&$0|1$&$2|0$&$2|0$&$2|0$\\
&$d_{4}$&=&$2\psa{1,1}3$&$1|1$&$3|1$&$3|1$&$3|1$\\
\hline
\end{tabular}
\end{center}
\caption{\protect Cohomology of $1|2$-Dimensional Complex Lie
Algebras}\label{Table 3}
\end{table}
\newpage
\section{The Moduli Space of $3|1$-dimensional Lie Superalgebras}

This is an interesting moduli space because it has a non nilpotent element, given by an extension of the simple Lie algebra $\sl(2,\C)$ by a $0|1$-dimensional
trivial algebra.  This is the algebra $d_1$ in the list of algebras below.  Other than this example, all $3|1$-dimensional Lie algebras are solvable, so
we next study how they are obtained by extensions.

\subsection{Construction of the moduli space of $3|1$-dimensional algebras}

Consider a $0|1$-dimensional vector space $W=\langle v_4\rangle$ and a $1|2$-dimensional vector space $M$. Note that the dimensions are reversed from the
algebra picture because we are studying the algebras as codifferentials on a $1|3$-dimensional space.   There are three nontrivial $1|2$-dimensional
codifferentials.  The generic $\lambda$ is of the form $\lambda=a_{11}\psa{1,4}1+a_{22}\psa{2,4}2+a_{32}\psa{2,4}3+a{2,3}\psa{3,4}2+a_{33}\psa{3,4}3$,
while $\psi$ must vanish.  The generic form of $\beta$ is $\beta=\pha{4}2b_1+\pha{4}3b_2$.

First we consider $\mu=2\psa{1,1}2+\psa{1,3}1-2\psa{2,3}2$.  The compatibility condition forces $a_{32}$ and $a_{33}$ to vanish and $a_{21}=-2a_{11}$. Moreover,
taking into account that we can add a term $[\mu,\beta]$ to $\lambda$, it turns out that we can assume that $\lambda=0$, so we obtain
the codifferential $d=\mu$, which is $d_{2}(1:0)$ on our list of algebras.

Next, consider the case $\mu=2\psa{11}2$. The compatibility condition forces $a_{32}=0$ and $a_{22}=-2a_{11}$. Moreover, $[\mu,\beta]=0$, so we are unable to eliminate any more terms in $\lambda$. The MC equation is automatically satisfied, so all of the possible forms of $\lambda$ above actually give codifferentials. The action of $G_{\delta,\mu}$ on $\lambda$ allows us to reduce to only three possible cases. The first is of the form
$\lambda=p\psa{1,4}1-2p\psa{2,4}2+\psa{3,4}2+q\psa{3,4}3$, $\lambda=\psa{1,4}1-2\psa{2,4}2-2\psa{3,4}3$, or $\lambda=0$.
The first choice of $\lambda$ gives the family $d_{2}(p:q)$, while the other 2 choices give $d_3$ and $d_4$, resp.

The last of the nonzero $1|2$-dimensional codifferentials is $\mu=p\psa{1,3}1+q\psa{2,3}2$.  The condition $[\mu,\lambda]=0$ has two possible solutions,
the first holding for generic values of $(p:q)$, while the second holds only for $(0:0)$.
In the first case, we have $a_{32}=a_{33}=0$, and when $q\ne 0$, we can eliminate the $a_{22}$ and $a_{23}$ terms by adding an appropriate $[\mu,\beta]$ term,
leaving only the coefficient $a_{11}$, which can be further reduced to the cases $a_{11}=1$ or $a_{11}=0$. When $a_{11}=1$, we obtain the codifferential
$d_5$, while when $a_{11}=0$, we obtain the subfamily $d_6(p:q:0)$ of the family $d_6(p:q:r)$. When $q=0$ and $p=1$, we substitute these values into $\mu$,
and now we will introduce new $p$ and $q$ variables by setting $a_{11}=p$ and $a_{22}=q$.  We still have a choice of $a_{23}$ which can be reduced to
being either 1 or 0.  As it turns out, the new variables $p$ and $q$ play a significant role, because when $a_{23}=1$ then if $q\ne0$ we obtain the codifferential
$d_5$, and if $q=0$ we obtain $d_{6}(1:0:0)$. Similarly, if $a_{23}=0$, then when $q\ne 0$, we obtain the codifferential $d_{5}$ and when $q=0$ we obtain
the codifferential $d_7(1: 0)$.

Returning to the original subcases, when $(p:q)=(0:0)$ in $\mu$, we obtain a more complicated solution because $\mu=0$.  Then $\lambda$ is given by a matrix
of the form $A=\left[\begin{array}{ccc}a_{11}&0&0\\0&a_{22}&a_{23}\\0&a_{32}&a_{33}\end{array}\right]$. Now, the compatibility and MC conditions are satisfied
automatically, so any $\lambda$ of this form gives a codifferential.  However, the group $G_{\delta,\mu}$ acts on the matrix of $\lambda$ by conjugating the
$2\times2$ submatrix of $\lambda$ by an element of $\gl(2,\C)$, and multiplying the matrix by a nonzero number.  This is a familiar pattern which says
that the submatrix can be reduced to a Jordan form,  in which case, we get two possibilities,
$A=\left[\begin{array}{ccc}p&0&0\\0&q&1\\0&0&r\end{array}\right]$ or $A=\left[\begin{array}{ccc}p&0&0\\0&q&0\\0&0&q\end{array}\right]$, the latter corresponding to an eigenvalue of geometric multiplicity 2 in the $2\times 2$ submatrix.  The first case gives $d_6(p:q:r)$, the second $d_7(p:q)$.

Finally, we have to consider the case $W=\langle v_1\rangle$ and $M=\langle v_2,v_3,v_4\rangle$, so that $W$ is a $1|0$-dimensional space and $M$ is a $0|3$-dimensional space.  There are 3 distinct elements of the moduli space of $0|3$-dimensional codifferentials, but the first one corresponds to the
simple Lie algebra, so we only have to consider the other two cases.  Now, for this space, $\lambda$ must vanish, but $\psi$ has a more complicated
form $$\psi=\psa{1,1}2c_1+\psa{1,1}3c_2+\psa{1,1}4c_3.$$

The first case is given by $\mu=\psa{2,4}2p+\psa{3,4}2+\psa{3,4}3q$.  The MC condition gives three solutions, the first for generic values of $p$ and $q$,
the second for $p=0$ and the third for $q=0$.   It is interesting to note that the $\mu$ is symmetric with respect to the interchange of $p$ and $q$, but
this symmetry is broken by the interaction with the $\psi$ term.  Generically, the codifferential is isomorphic to $d_6(0:p:q)$, while if $p=0$, then the
codifferential is isomorphic to $d_2(0:q)$ if $c_1\ne 0$ or $d_6(0:0:q)$, when $c_1=0$. On the other hand, if $q=0$, then the codifferential
is isomorphic to $d_2(0:p)$ when $c_1\ne 0$ and $d_6(0: 0:p)$ when $c_1=0$. Note that the symmetry in the isomorphism classes between $p$ and $q$ is
restored in the isomorphism class, which has to be the case because $\mu(p:q)\sim\mu(q:p)$, so extensions should not be different when $p$ vanishes
than when $q$ vanishes.

The second case is given by $\mu=\psa{2,4}2+\psa{3,4}3$. The MC equation forces $\psi=0$, so $d=\mu$, which is isomorphic to $d_7(p:q)$.

Finally, the last case is given by $\mu=0$, in which case the codifferential $d$ is just $\psi$ and there are no conditions.  On the other
hand, it is easy to see that whenever any of the terms in $\psi$ is nonzero, we obtain an equivalent codifferential, which is isomorphic to $d_4$.
Otherwise, we just obtain the zero codifferential, which is also given by $d_7(0:0)$.

In the table of algebras below, we include some special subfamilies of $d_6(p:q:r)$. Because this algebra is parametrized by $\P^2/\Sigma_2$, where
the action of $\Sigma_2$ on $\P^2$ is given by interchanging the second two coordinates, there are special $\P^1$s for which the cohomology
does not follow the generic pattern.  There are also some special points in the families for which the generic pattern does not hold. We include
the special points of the $d_2(p:q)$ family in the main table, but list the special points for $d_6(p:q:r)$ and $d_7(p:q)$ in separate tables.

\begin{table}[ht]
\begin{center}
\begin{tabular}{rlllrrrr}
&Algebra&&Codifferential&$h_0$&$h_1$&$h_2$&$h_3$\\
\hline\\
&$d_{1}$&=&$\psa{2,3}4+\psa{2,4}3+\psa{3,4}2$&$1|0$&$1|0$&$1|0$&$1|1$\\
&$d_{2}(p:q)$&=&$8\psa{1,1}2+p\psa{1,4}1-2p\psa{2,4}2+\psa{3,4}2+q\psa{3,4}3$&$0|0$&$1|0$&$0|1$&$0|0$\\
&$d_{2}(0:1)$&=&$8\psa{1,1}2+\psa{3,4}2+\psa{3,4}3$&$0|1$&$2|0$&$0|1$&$0|0$\\
&$d_{2}(1:1)$&=&$8\psa{1,1}2+\psa{1,4}1-2\psa{2,4}2+\psa{3,4}2+\psa{3,4}3$&$0|0$&$1|0$&$1|1$&$0|1$\\
&$d_{2}(0:0)$&=&$8\psa{1,1}2+\psa{3,4}2$&$0|1$&$4|2$&$4|4$&$1|2$\\
&$d_{3}$&=&$8\psa{1,1}2+\psa{1,4}1-2\psa{2,4}2-2\psa{3,4}3$&$0|0$&$2|0$&$0|2$&$0|0$\\
&$d_{4}$&=&$8\psa{1,1}2$&$0|3$&$7|2$&$4|5$&$1|2$\\
&$d_{5}$&=&$\psa{2,3}2+\psa{1,4}1$&$0|0$&$0|0$&$0|0$&$0|0$\\
&$d_{6}(p:q:r)$&=&$p\psa{1,4}1+q\psa{2,4}2+\psa{3,4}2+r\psa{3,4}3$&$0|0$&$2|0$&$0|2$&$0|0$\\
&$d_{6}(p:q:p+q)$&=&$p\psa{1,4}1+q\psa{2,4}2+\psa{3,4}2+(p+q)\psa{3,4}3$&$0|0$&$2|0$&$1|2$&$0|1$\\
&$d_{6}(p:q:-p-q)$&=&$p\psa{1,4}1+q\psa{2,4}2+\psa{3,4}2-(p+q)\psa{3,4}3$&$0|0$&$2|0$&$1|2$&$0|1$\\
&$d_{6}(p:q:0)$&=&$p\psa{1,4}1+q\psa{2,4}2+\psa{3,4}2$&$0|1$&$3|0$&$0|3$&$1|0$\\
&$d_{6}(p:q:-2p)$&=&$p\psa{1,4}1+q\psa{2,4}2+\psa{3,4}2-2p\psa{3,4}3$&$0|0$&$2|0$&$0|3$&$1|0$\\
&$d_{6}(0:p:q)$&=&$p\psa{2,4}2+\psa{3,4}2+q\psa{3,4}3$&$1|0$&$2|1$&$2|2$&$2|2$\\
&$d_{6}(p:p:q)$&=&$p\psa{1,4}1+p\psa{2,4}2+\psa{3,4}2+q\psa{3,4}3$&$0|0$&$2|0$&$0|2$&$0|1$\\
&$d_{6}(p:q:2p+q)$&=&$p\psa{1,4}1+q\psa{2,4}2+\psa{3,4}2+(2p+q)\psa{3,4}3$&$0|0$&$2|0$&$0|2$&$1|0$\\
&$d_{6}(p:q:-3p)$&=&$p\psa{1,4}1+q\psa{2,4}2+\psa{3,4}2-3p\psa{3,4}3$&$0|0$&$2|0$&$0|2$&$0|1$\\
&$d_{7}(p:q)$&=&$p\psa{1,4}1+q\psa{2,4}2+q\psa{3,4}3$&$0|0$&$4|0$&$0|4$&$0|0$\\
\hline
\end{tabular}
\end{center}
\caption{\protect Cohomology of $3|1$-Dimensional Complex Lie
Algebras}\label{Table 4}
\end{table}

\subsection{Deformations of the $3|1$-dimensional algebras}

The algebra $d_1$ is rigid.  The family $d_2(p:q)$ generically only deforms along the family.  There are two special points $(1:1)$ and $(0:1)$ where the
dimension of $H^1$ is 2, rather than the generic value 1, but this does not affect the deformations.  On the other hand, $d_2(0:0)$ has jump deformations to
all the other elements in the family $d_2(p:q)$, so the dimension of its $H^2$ is not generic.  The algebra $d_3$ has a jump deformation to $d_2(1:-2)$ and
also deforms smoothly in a nbd of $d_2(1:-2)$. In fact, if there is a jump deformation from an algebra to a member of a family, then there always smooth
deformations in a nbd of this point. The algebra $d_4$ has jump deformations to $d_2(x:y)$ for all $(x:y)$ as well as a jump to $d_3$. The algebra $d_5$ is
completely rigid.

\begin{table}[ht]
\begin{center}
\begin{tabular}{rlllrrrr}
&Algebra&&Codifferential&$h_0$&$h_1$&$h_2$&$h_3$\\
\hline\\
&$d_{6}(1:-1:0)$&=&$\psa{1,4}1-\psa{2,4}2+\psa{3,4}2$&$0|1$&$3|2$&$4|3$&$1|2$\\
&$d_{6}(1:-2:0)$&=&$\psa{1,4}1-2\psa{2,4}2+\psa{3,4}2$&$0|1$&$3|0$&$0|4$&$3|0$\\
&$d_{6}(1:-2:-1)$&=&$\psa{1,4}1-2\psa{2,4}2+\psa{3,4}2-\psa{3,4}3$&$0|0$&$2|2$&$3|3$&$1|1$\\
&$d_{6}(1:-3:-2)$&=&$\psa{1,4}1-3\psa{2,4}2+\psa{3,4}2-2\psa{3,4}3$&$0|0$&$2|0$&$1|3$&$1|2$\\
&$d_{6}(1:1:3)$&=&$\psa{1,4}1+\psa{2,4}2+\psa{3,4}2+3\psa{3,4}3$&$0|0$&$2|0$&$0|2$&$1|1$\\
&$d_{6}(2:-1:2)$&=&$2\psa{1,4}1-\psa{2,4}2+\psa{3,4}2+2\psa{3,4}3$&$0|0$&$2|0$&$0|2$&$0|1$\\
&$d_{6}(0:0:1)$&=&$\psa{3,4}2+\psa{3,4}3$&$1|1$&$3|3$&$4|4$&$4|4$\\
&$d_{6}(0:1:-1)$&=&$\psa{2,4}2+\psa{3,4}2-\psa{3,4}3$&$1|0$&$2|1$&$3|3$&$4|4$\\
&$d_{6}(1:0:1)$&=&$\psa{1,4}1+\psa{3,4}2+\psa{3,4}3$&$0|1$&$3|0$&$1|3$&$1|2$\\
&$d_{6}(1:-2:1)$&=&$\psa{1,4}1-2\psa{2,4}2+\psa{3,4}2+\psa{3,4}3$&$0|0$&$2|0$&$1|3$&$1|2$\\
&$d_{6}(1:-1:1)$&=&$\psa{1,4}1-\psa{2,4}2+\psa{3,4}2+\psa{3,4}3$&$0|0$&$2|2$&$2|2$&$1|1$\\
&$d_{6}(1:-3:1)$&=&$\psa{1,4}1-3\psa{2,4}2+\psa{3,4}2+\psa{3,4}3$&$0|0$&$2|0$&$0|2$&$0|2$\\
&$d_{6}(1:-5:-3)$&=&$\psa{1,4}1-5\psa{2,4}2+\psa{3,4}2-3\psa{3,4}3$&$0|0$&$2|0$&$0|2$&$1|1$\\
&$d_{6}(0:1:1)$&=&$\psa{2,4}2+\psa{3,4}2+\psa{3,4}3$&$1|0$&$2|1$&$2|2$&$2|2$\\
&$d_{6}(1:1:2)$&=&$\psa{1,4}1+\psa{2,4}2+\psa{3,4}2+2\psa{3,4}3$&$0|0$&$2|0$&$1|2$&$0|2$\\
&$d_{6}(1:-4:-3)$&=&$\psa{1,4}1-4\psa{2,4}2+\psa{3,4}2-3\psa{3,4}3$&$0|0$&$2|0$&$1|2$&$0|2$\\
&$d_{6}(2:-3:1)$&=&$2\psa{1,4}1-3\psa{2,4}2+\psa{3,4}2+\psa{3,4}3$&$0|0$&$2|0$&$1|2$&$1|1$\\
&$d_{6}(1:-3:2)$&=&$\psa{1,4}1-3\psa{2,4}2+\psa{3,4}2+2\psa{3,4}3$&$0|0$&$2|0$&$1|2$&$0|2$\\
&$d_{6}(1:-4:-2)$&=&$\psa{1,4}1-4\psa{2,4}2+\psa{3,4}2-2\psa{3,4}3$&$0|0$&$2|0$&$0|3$&$2|0$\\
&$d_{6}(1:0:2)$&=&$\psa{1,4}1+\psa{3,4}2+2\psa{3,4}3$&$0|1$&$3|0$&$0|3$&$2|0$\\
&$d_{6}(1:-3:0)$&=&$\psa{1,4}1-3\psa{2,4}2+\psa{3,4}2$&$0|1$&$3|0$&$0|3$&$1|1$\\
&$d_{6}(1:0:0)$&=&$\psa{1,4}1+\psa{3,4}2$&$0|1$&$3|0$&$0|3$&$1|0$\\
&$d_{6}(0:0:0)$&=&$\psa{3,4}2$&$1|1$&$5|3$&$7|8$&$9|9$\\
\hline
\end{tabular}
\end{center}
\caption{\protect Cohomology of special points in the family $d_6(p:q:r)$}\label{Table 5}
\end{table}

The family of algebras $d_6(p:q:r)$ generically has $h_2=2$, which is precisely what one would expect as it deforms only along the family. For the special
subfamilies, most of them deform generically, but $d_6(p:q:0)$ also has a jump deformation to $d_5$, and $d_6(p:q:-2p)$ also has a jump to $d_2(p:q)$ as
well as deforming in a nbd of $d_2(p:q)$. Now, for the special points, if they belong to a special subfamily, then they take part in any extra deformations
which the subfamily has. However, one has to keep in mind that the symmetry of interchanging the last two coordinates must be taken into account, as was
done when we listed the subfamilies.  For example, the subfamily $d_6(p:q:0)$ coincides (up to isomorphism) with the subfamily $d_6(p:0:q)$, so we only listed
the first subfamily.  Similarly, the element $d_6(1:-2:0)$ is equivalent to $d_6(1:0: -2)$, so it both jumps to $d_5$ and to $d_2(1:0)$.  Other than
these cases, there is only one additional special case, the element $d_6(0:1:-1)$, which has a jump deformation to $d_1$. Finally, the generic element $d_6(0:0:0)$ has jump deformations to $d_1$, $d_2(x: y)$ for all $(x:y)$, $d_5$ and $d_6(x:y:z)$ except $(0:0:0)$.  Another important observation is that
if any element of a family has a deformation to some algebra, then the generic element also has such a deformation.  This is a case of the observation that
if algebra $a$ deforms to algebra $b$ and algebra $b$ deforms to algebra $c$, then algebra $a$ also deforms to algebra $c$.

\begin{table}[ht]
\begin{center}
\begin{tabular}{rlllrrrr}
&Algebra&&Codifferential&$h_0$&$h_1$&$h_2$&$h_3$\\
\hline\\
&$d_{7}(1:0)$&=&$\psa{1,4}1$&$0|2$&$6|0$&$0|6$&$2|0$\\
&$d_{7}(2:-1)$&=&$2\psa{1,4}1-\psa{2,4}2-\psa{3,4}3$&$0|0$&$4|0$&$1|4$&$0|1$\\
&$d_{7}(1:-2)$&=&$\psa{1,4}1-2\psa{2,4}2-2\psa{3,4}3$&$0|0$&$4|0$&$0|6$&$2|0$\\
&$d_{7}(0:1)$&=&$\psa{2,4}2+\psa{3,4}3$&$1|0$&$4|1$&$4|4$&$4|4$\\
&$d_{7}(1:-3)$&=&$\psa{1,4}1-3\psa{2,4}2-3\psa{3,4}3$&$0|0$&$4|0$&$0|4$&$0|2$\\
&$d_{7}(1:1)$&=&$\psa{1,4}1+\psa{2,4}2+\psa{3,4}3$&$0|0$&$4|0$&$0|4$&$0|2$\\
&$d_{7}(1:-1)$&=&$\psa{1,4}1-\psa{2,4}2-\psa{3,4}3$&$0|0$&$4|4$&$4|4$&$0|0$\\
&$d_{7}(0:0)$&=&$0$&$1|3$&$10|6$&$13|15$&$16|16$\\
\hline
\end{tabular}
\end{center}
\caption{\protect Cohomology of special points in the family $d_7(p:q)$}\label{Table 6}
\end{table}

Generically, the algebra $d_7(p:q)$ has jump deformations to $d_6(p:q:q)$ and deforms in a nbd of $d_6(p:q:q)$ and $d_7(p:q)$. The special point
$d_7(1:0)$ also has a jump deformation to $d_5$, while the special point $d_7(1:-2)$ has jump deformations to $d_2(x:y)$ for all $(x:y)$ except

as to $d_3$. This family is parametrized by $\P^1$, with no action of a symmetric group. The generic element $d_7(0:0)$ is just the trivial algebra,
so it has jump deformations to every other element in the moduli space.

\section{The Moduli Space of $2|2$-dimensional Lie Superalgebras}
\subsection{Construction of the moduli space of $2|2$-dimensional algebras}
This is the most complicated of the moduli spaces in this paper to construct.  We begin by considering extensions of the trivial algebra
structure $\delta=0$ on the $0|1$-dimensional space $W=\langle v_4\rangle$ by an algebra structure $\mu$ on a $2|1$-dimensional space.  The
generic $\lambda$ is of the form
$$\lambda=\psa{1,4}1a_{11}+\psa{2,4}1a_{12}+\psa{1,4}2a_{21}+\psa{2,4}2a_{22}+\psa{3,4}3a_{33},$$ which has a block diagonal matrix
$A=\left[\begin{array}{ccc}a_{11}&a_{12}&0\\a_{21}&a_{22}&0\\0&0&a_{33}\end{array}\right]$.  The $\psi$ term must vanish, and $\beta=\pha{4}3b$.

There are 4 nontrivial possibilities for $\mu$. The first case is given by $\mu=\psa{1,3}1p+\psa{2,3}1+\psa{2,3}2q$.  There are two solutions to the
compatibility condition, the first holding for generic values of $p$ and $q$, and the second holding only for $p=q=0$. For the first solution,
we have $a_{21}=a_{31}=0$ and $a_{11}=pa_{12}+a_{22}-a_{12}q$, so only two free variables remain, which we will denote by $r=a_{22}$ and $s=a_{12}$ for short.
The MC condition is satisfied automatically, so $d=\mu+\lambda$.  When $p\ne q$ and $r\ne sq$, this gives the algebra $d_1$, when $p=q$ and $r\ne sq$, this
gives the algebra $d_5(p:0)$, and finally when $r=sq$, we obtain the algebra $d_10(p:q:0)$.

Next, consider the case when both $p$ and $q$ vanish.  Then the compatibility condition yields $a_{21}=0$, and adding a coboundary term allows us to eliminate
the coefficient $a_{12}$ as well, leaving two coefficients $a_{11}=p$ and $a_{22}=q$. The reader may be curious why we make this $p$ and $q$ substitution,
and the explanation is that considering the action of $G_{\delta,\mu}$ on $\lambda$, we note that the two coefficients are scaled by a number, which suggests
that they represent projective coefficients. Not all such projectively given coefficient relations survive the test of isomorphism, but in this case, they do,
as the algebra we have constructed is isomorphic to $d_5(p:q)$, except that when $(p:q)=(0:0)$, we discover that $d_5(0:0)\sim d_{10}(0:0:0)$, a type of
occurrence which often happens in our construction of moduli spaces of algebras.

Next, consider $\mu=\psa{1,3}1+\psa{2,3}2$. The compatibility condition forces $a_{33}=0$ and by adding a coboundary, we could eliminate either the
$a_{11}$ or the $a_{22}$ term, but this time, we found it convenient not to use this simplification, for reasons we will explain. The matrix of $\lambda$
is essentially a $2\times 2$ block matrix, and the action of $G_{\delta,\mu}$ on $\lambda$ is essentially given by conjugation of the $2\times 2$ submatrix
of $\lambda$, up to a constant multiple. This action we know well, and it gives isomorphism classes of Jordan decompositions, so we know the decomposition
can be reduced to some simple cases.  The first is when the submatrix is of the form $A=\left[\begin{array}{cc}p&1\\0&q\end{array}\right]$. In this case,
whenever $p\ne q$, we obtain $d_1$, and when $p=q$, we obtain $d_{5}(p:0)$.  The second is given by the diagonal matrix $\diag(p,p)$, and in this case,
independently of $p$, we obtain $d_{11}(1:0)$.

Now, let $\mu=\psa{1,2}3+2\psa{1,1}3$. The compatibility condition forces $a_{12}=0$ and $a_{11}= -a_{22}-a_{33}$, but there is no additional
simplification by adding a coboundary term, since $[\mu,\beta]$ vanishes. The matrix of $\lambda$ is determined by two coefficients. This time
we found it convenient to set $a_{22}=p$ and $a_{33}=-p-q$. The MC condition is automatically satisfied, and we obtain that the algebra is isomorphic
to $d_6(p:q)$.

The last nontrivial $\mu$ is $\mu=\psa{1,1}3$. The compatibility conditions force $a_{12}=0$ and $a_{33}=-a_{11}$. Adding a coboundary term does not
change anything, but applying $G_{\delta,\mu}$, we discover that unless $a_{11}=a_{22}$, we can eliminate the $a_{21}$ term, and in this case,
we obtain $d_{11}(p:q)$, except when $p=q$. When $p=q\ne0$, we obtain $d_8$, and when $p=q=0$, we obtain $d_9$.  When $a_{11}=a_{22}$, and $a_{21}\ne0$,
then when $a_{21}\ne 0$, we obtain $d_7(a_{11},a_{11})$, and otherwise if $a_{11}\ne0$, we get $d_8$ and when $a_{11}=0$ we obtain $d_9$.

Finally, we need to analyze the case $\mu=0$. Then $\lambda$ is given by the generic value of $\lambda$, whose matrix is the block diagonal $A$ given
in the beginning of this section, and the group $G_{\delta,\mu}$ acts by reducing to Jordan form.  Thus there are two cases to consider, given by
the matrices below.
$$\left[\begin{array}{ccc}p&1&0\\0&q&0\\0&0&r\end{array}\right],\qquad \left[\begin{array}{ccc}p&0&0\\0&p&0\\0&0&q\end{array}\right].$$
The first matrix gives the algebra $d_{10}(p:q:r)$, while the second gives $d_{11}(p:q)$.

Next, we extensions of the trivial algebra
structure $\delta=0$ on the $1|0$-dimensional space $W=\langle v_2\rangle$ by an algebra structure $\mu$ on a $1|2$-dimensional space
$M=\langle v_1,v_3,v_4\rangle$.  The
generic $\lambda$ is of the form $$\lambda=\psa{1,2}3a_{21}+\psa{1,2}4a_{31}+\psa{2,3}1a_{12}+\psa{2,4}1a_{13}.$$  Generically, $\psi=\psa{2,2}3c_1+\psa{2,2}4c_2$
and $\beta=\pha{2}1b$. There are three nontrivial possibilities for $\mu$.

The first case is $\mu=4\psa{1,1}3+\psa{1,4}1-2\psa{2,3}1$.  The compatibility condition gives $a_{12}=a_{31}=0$ and $a_{21}=8a_{13}$, but adding a coboundary
term allows us to eliminate $\lambda$. The MC condition forces $\psi=0$, so the resulting algebra is just $\mu$, which is isomorphic to $d_7(1:0)$.

The second case is $\mu=4\psa{1,1}3$. The compatibility condition gives $a_{12}=a_{13}=0$, and adding a coboundary term allows us to assume $a_{12}=0$,
leaving only the coefficient $a_{31}$ in $\lambda$, which can be taken to be 1 or 0.  The MC condition and the cocycle condition
are automatically satisfied, and thus we have no conditions on
the coefficients $c_1$ or $c_2$.  Assuming $a_{31}=1$, then when $c_1=-4c_2^2$, the algebra is isomorphic to $d_2$, and otherwise, it is isomorphic to $d_{4}$.
When $a_{31}=0$, if $c_2\ne 0$, it is $d_2$,  if $c_2=0$ and $c_1\ne0$, then we obtain $d_6(0:0)$, while if $c_1=c_2=0$, we obtain $d_9$.

Finally, when $\mu=\psa{1,4}1p+\psa{3,4}3q$ (note this includes the case $\mu=0$), the compatibility condition gives three solutions, a generic one, one for $q=-p$, and one where $p=q=0$.  Let us consider the generic case first, which forces $a_{12}=a_{21}=a_{31}=0$.

The MC condition gives three solutions, depending on $p$ and $q$.  The first case is generic in terms of $p$ and $q$, and $\psi$ must vanish. The resulting
algebra is isomorphic to $d_{10}(p:0:q)$ except that if $p=0$, we must have $a_{1,3,1}\ne 0$. When $p=a_{1,3}=0$, we obtain $d_{11}(0:q)$, in other words,
we get $d_{11}(0:1)$ and $d_{11}(0:0)$.

The second solution to the MC condition
 has $q=0$ and $c_2=0$.  This one gives rise  to several different algebras, which is not surprising, since we have the variables $c_1$, $a{1,3,1}$ and $p$ to consider. Depending on the values of these coefficients, we get $d_7(0:1)$, $d_7(0:0)$, $d_9$, $d_{10}(1:0:0)$, $d_{10}(0:0:0)$ and $d_{11}(0:0)$.

The final solution to the MC condition occurs when both $p$ and $q$ vanish, and this time, the cocycle condition forces either $c_2$ or $a_{1,3}$ to vanish.
In the first case, we have two potentially nonzero coeffients, $c_1$ and $a_{1,3}$, and four algebras arise, $d_7(0:0)$, $d_9$, $d_{10}(0:0:0)$, and $d_{11}(0:0)$.

The second solution of the compatibility condition has $q=-p$ and $a_{31}=0$, leaving $a_{21}$, $a_{12}$ and $a_{1,3}$ undetermined. When $p\ne0$ we can assume it is 1, and in  this case, we can eliminate the coefficient $a_{13}$ by
adding a coboundary term. By applying an element of $G_{\delta,\mu}$, we can also reduce to the cases where $a_{21}$ and $a_{12}$ are either 1 or 0.
When they both equal 1, the MC condition forces $c_1=0$ and $c_2=-1$, completely determining the algebra, which is $d_{3}$. When $a_{21}=0$ and $a_{12}=1$,
the MC condition forces $c_1=c_2=0$, giving the algebra $d_5(0:1)$. When $a_{21}=1$ and $a_{12}=0$, the MC condition again forces $c_1=c_2=0$, giving the
algebra $d_6(1:0)$. Finally, when $a_{21}=a_{12}=0$, the MC condition also forces $c_1=c_2=0$, giving the algebra $d_{10}(1:0:-1)$.

Now, we study the case when $p=0$. The MC condition has two solutions, $a_{21}=0$ or both $a_{12}$ and $a_{13}$ vanish.  When $a_{21}=0$, the cocycle condition
gives $c_1a_{12}+c_2a_{13}=0$. We can break the solution to the cocycle condition into 3 parts. After some analysis, we
obtain the algebras $d_6$, $d_7(0:0)$, $d_9$, $d_{10}(0:0:0)$, and $d_{11}(0:0)$. The other solution to the MC condition givess $a_{12}=a_{13}=0$, and the
cocycle condition puts no restriction on $\psi$. Depending on the values of $a_{21}$, $c_1$ and $c_2$, we obtain the algebras $d_4$, $d_6$, $d_9$ and $d_{11}(0:0)$.

\begin{table}[ht]
\begin{center}
\begin{tabular}{rlllrrrr}
&Algebra&&Codifferential&$h_0$&$h_1$&$h_2$&$h_3$\\
\hline\\
&$d_{1}$&=&$\psa{1,3}1+\psa{2,3}1+2\psa{2,3}2+\psa{2,4}1+\psa{2,4}2$&$0|0$&$0|0$&$0|0$&$0|0$\\
&$d_{2}$&=&$\psa{1,2}3+4\psa{1,1}3+\psa{1,2}4+2\psa{2,2}4$&$0|2$&$2|2$&$2|0$&$0|0$\\
&$d_{3}$&=&$\psa{1,2}3+4\psa{1,1}3+2\psa{1,4}2+\psa{2,4}2$
\\&&&$-\psa{3,4}3+\psa{1,2}4-\psa{2,3}2-\psa{3,4}4$&$0|0$&$0|1$&$1|0$&$0|0$\\
&$d_{4}$&=&$8\psa{1,1}3+\psa{1,2}4+2\psa{1,1}4$&$0|2$&$3|2$&$3|2$&$2|2$\\
&$d_{5}(p:q)$&=&$\psa{2,3}1+p-q\psa{1,4}1+p\psa{2,4}2+q\psa{3,4}3$&$0|0$&$1|0$&$0|1$&$0|0$\\
&$d_{6}(p:q)$&=&$\psa{1,2}3+4\psa{1,1}3+q\psa{1,4}1+2(p-q)\psa{1,4}2$
\\&&&$+p\psa{2,4}2-(p+q)\psa{3,4}3$&$0|0$&$1|0$&$0|1$&$0|0$\\
&$d_{7}(p:q)$&=&$4\psa{1,1}3+p\psa{1,4}1+\psa{1,4}2+q\psa{2,4}2-2p\psa{3,4}3$&$0|0$&$1|0$&$0|1$&$0|0$\\
&$d_{8}$&=&$4\psa{1,1}3+\psa{1,4}1+\psa{2,4}2-2\psa{3,4}3$&$0|0$&$2|0$&$0|3$&$1|0$\\
&$d_{9}$&=&$4\psa{1,1}3$&$1|2$&$5|4$&$6|6$&$6|6$\\
&$d_{10}(p:q:r)$&=&$p\psa{1,4}1+\psa{2,4}1+q\psa{2,4}2+r\psa{3,4}3$&$0|0$&$2|0$&$0|2$&$0|0$\\
&$d_{11}(p:q)$&=&$p\psa{1,4}1+p\psa{2,4}2+q\psa{3,4}3$&$0|0$&$4|0$&$0|4$&$0|0$\\
\hline
\end{tabular}
\end{center}
\caption{\protect Cohomology of $2|2$-Dimensional Complex Lie
Algebras}\label{Table 7}
\end{table}
The third and last solution of the compatibility condition gives $p=q=0$, in other words, $\mu=0$. There are 2 solutions to the MC condition,
$a_{21}=a_{31}=0$ or $a_{12}=a_{13}=0$. In the first case, the cocycle condition gives $c_1a_{12}+c_2a_{13}=0$, which can be divided into 3 cases,
when $c_1\ne0$, when $c_1=c_2=0$ and when $c_1=a_{13}=0$. The first gives $d_7(0:0)$ when $a_{13}\ne0$ and $d_9$ when $a_{13}=0$. The second gives
$d_{10}(0:0:0)$ as long as not both $a_{12}$ and $a_{13}$ vanish, and $d_{11}(0:0)$ otherwise.  The third solution gives $d_7(0:0)$ when neither $a_{12}$ nor
$c_2$ vanish, $d_{10}(0:0:0)$ when $a_{12}\ne0$ and $c_2=0$, $d_9$ when $a_{12}=0$ and $c_2\ne0$ and $d_{11}(0:0)$ when both $a_{12}$ and $c_2$ vanish.
For the second solution to the MC condition, the cocycle condition is trivial. The solutions, depending on the values of the coefficients $a{21}$, $a_{31}$,
$c_1$ and $c_2$, are $d_4$, $d_6(0:0)$, $d_9$, and  $d_{11}(0:0)$.

This completes the construction of the moduli space of $2|2$-dimensional algebras.

\subsection{Deformations of the $2|2$-dimensional algebras}
The algebras $d_1$, $d_2$ and $d_3$ are rigid, although only $d_1$ is totally rigid (cohomology vanishes identically).
The algebra $d_4$ has jump deformations to both $d_2$ and $d_3$.

The family $d_5(p:q)$ is parameterized by $\P^1$, with no action of a symmetric group. Generically, deformations are
only along the family. The special point $d_5(1:0)$ has an
additional jump deformation to $d_3$, while the special point $d_5(0:1)$ has a jump deformation to $d_1$.
The points $d_5(1:1)$, $d_5(2:1)$ and $d_5(3:1)$ have cohomology dimensions that do not fit the generic
picture, but deform generically.  The generic point $d_5(0:0)$ is quite special. It is isomorphic to
$d_{10}(0:0:0)$, and it has jump deformations to $d_1$, $d_3$, $d_5(x: y)$ except $(0: 0)$,
$d_6(x: y)$ except $(1: 1)$ and $(0: 0)$, $d_7(x: y)$ and $d_{10}(x: y: z)$ except $(0: 0: 0)$.

\begin{table}[ht]
\begin{center}
\begin{tabular}{rlllrrrr}
&Algebra&&Codifferential&$h_0$&$h_1$&$h_2$&$h_3$\\
\hline\\
&$d_{5}(p:q)$&=&$\psa{2,3}1+p-q\psa{1,4}1+p\psa{2,4}2+q\psa{3,4}3$&$0|0$&$1|0$&$0|1$&$0|0$\\
&$d_{5}(1:0)$&=&$\psa{2,3}1+\psa{1,4}1+\psa{2,4}2$&$0|0$&$1|0$&$0|2$&$1|0$\\
&$d_{5}(0:1)$&=&$\psa{2,3}1-\psa{1,4}1+\psa{3,4}3$&$0|0$&$1|1$&$2|2$&$2|2$\\
&$d_{5}(1:1)$&=&$\psa{2,3}1+\psa{2,4}2+\psa{3,4}3$&$1|0$&$1|1$&$1|1$&$1|1$\\
&$d_{5}(3:1)$&=&$\psa{2,3}1+2\psa{1,4}1+3\psa{2,4}2+\psa{3,4}3$&$0|0$&$1|0$&$0|1$&$0|1$\\
&$d_{5}(2:1)$&=&$\psa{2,3}1+\psa{1,4}1+2\psa{2,4}2+\psa{3,4}3$&$0|0$&$1|0$&$0|1$&$0|0$\\
&$d_{5}(0:0)$&=&$\psa{2,3}1$&$1|1$&$4|3$&$5|6$&$6|6$\\
\hline\\
\end{tabular}
\end{center}
\caption{\protect Cohomology of the family $d_5(p:q)$}
\label{Table X}
\end{table}

The family $d_6(p:q)$ is parametrized by $\P^1/\Sigma_2$, where $\Sigma_2$ acts by permuting the coordinates $(p:q)$.
Generically, an point in this family only has deformations in a nbd of the point. The special point $d_6(1:0)$ also
has a jump deformation to $d_3$, while the points $d_6(1:1)$ and $d_6(1:-1)$ deform generically.  The generic point
$d_6(0:0)$ has jump deformations to $d_2$, $d_3$, $d_4$, and $d_6(x: y)$ except $(0: 0)$. The cohomology of the elements
in this family is given in the table below.

\begin{table}[ht]
\begin{center}
\begin{tabular}{rlllrrrr}
&Algebra&&Codifferential&$h_0$&$h_1$&$h_2$&$h_3$\\
\hline\\
&$d_{6}(p:q)$&=&$\psa{1,2}3+4\psa{1,1}3+q\psa{1,4}1$\\
&&&$+2(p-q)\psa{1,4}2+p\psa{2,4}2-(p+q)\psa{3,4}3$&$0|0$&$1|0$&$0|1$&$0|0$\\
&$d_{6}(1:0)$&=&$\psa{1,2}3+4\psa{1,1}3+2\psa{1,4}2+\psa{2,4}2-\psa{3,4}3$&$0|0$&$1|1$&$2|2$&$2|2$\\
&$d_{6}(1:-1)$&=&$\psa{1,2}3+4\psa{1,1}3-\psa{1,4}1+4\psa{1,4}2+\psa{2,4}2$&$0|1$&$2|0$&$0|1$&$0|0$\\
&$d_{6}(0:0)$&=&$\psa{1,2}3+4\psa{1,1}3$&$0|2$&$4|2$&$4|4$&$4|4$\\
\hline\\
\end{tabular}
\end{center}
\caption{\protect Cohomology of elements in the family $d_6(p:q)$}
\label{Table Y}
\end{table}

The family $d_7(p:q)$ is parametrized by $\P^1$, with no action of a symmetric group. Generically, elements in
the family deform only along the family. The special point $d_7(1:1)$ deforms in a nbd of $d_6(1:1)$ but does
 not have a jump deformation to $d_6(1:1)$. The special points $d_7(1:0)$, $d_7(0:1)$, $d_7(1:-1)$ and
$d_7(3:2)$ deform generically.  The generic element $d_7(0:0)$ has jump deformations to
$d_3$, $d_6(x: y)$ except $(0: 0)$, and $d_7(x: y)$ except $(0: 0)$.

\begin{table}[ht]
\begin{center}
\begin{tabular}{rlllrrrr}
&Algebra&&Codifferential&$h_0$&$h_1$&$h_2$&$h_3$\\
\hline\\
&$d_{7}(p:q)$&=&$4\psa{1,1}3+p\psa{1,4}1+\psa{1,4}2+q\psa{2,4}2-2p\psa{3,4}3$&$0|0$&$1|0$&$0|1$&$0|0$\\
&$d_{7}(1:0)$&=&$4\psa{1,1}3+\psa{1,4}1+\psa{1,4}2-2\psa{3,4}3$&$1|0$&$1|1$&$1|1$&$1|1$\\
&$d_{7}(0:1)$&=&$4\psa{1,1}3+\psa{1,4}2+\psa{2,4}2$&$0|1$&$2|0$&$0|1$&$0|0$\\
&$d_{7}(1:1)$&=&$4\psa{1,1}3+\psa{1,4}1+\psa{1,4}2+\psa{2,4}2-2\psa{3,4}3$&$0|0$&$1|0$&$0|2$&$1|0$\\
&$d_{7}(1:-1)$&=&$4\psa{1,1}3+\psa{1,4}1+\psa{1,4}2-\psa{2,4}2-2\psa{3,4}3$&$0|0$&$1|0$&$0|1$&$0|0$\\
&$d_{7}(3:2)$&=&$4\psa{1,1}3+3\psa{1,4}1+\psa{1,4}2+2\psa{2,4}2-6\psa{3,4}3$&$0|0$&$1|0$&$0|1$&$0|1$\\
&$d_{7}(0:0)$&=&$4\psa{1,1}3+\psa{1,4}2$&$1|1$&$3|3$&$4|4$&$4|4$\\
\hline\\
\end{tabular}
\end{center}
\caption{\protect Cohomology of the elements of the family $d_7(p:q)$}
\label{Table Z}
\end{table}

The algebra $d_8$ has jump deformations to $d_6(1:1)$, $d_7(1:1)$ and deforms in a nbd of each of these points.

The algebra $d_9$ has jump deformations to $d_3$, $d_4$, and $d_6(x: y)$ and $d_7(x: y)$ for all $(x:y)$ as well as $d_8$.

The family $d_{10}(p:q:r)$ is parametrized by $\P^2/\Sigma_2$, where $\Sigma_2$ acts by permuting the first two coordinates.
Generically, elements deform only along the family. There are some special subfamilies parametrized by $\P^1$, for which the
cohomology or deformation theory is not generic, and also some special points, each of which belongs to one or more of the
special subfamilies, for which the cohomology or deformation theory is even more unusual.

The members of the subfamily $d_{10}(p:q:p-q)$ have deformations in a nbd of $d_5(p:p-q)$, but don't have jump deformations to this point.
The members of the subfamily $d_{10}(p:q:0)$ all have jump deformations to $d_1$.  The members of the subfamily $d_{10}(p:q:-2p)$ have
jump deformations to $d_7(p:q)$ and deform in a nbd of this point, while the members of the subfamily $d_{10}(p:q:-p-q)$ jump to $d_6(p:q)$
and deform in a nbd of this point. The other special subfamilies have non generic cohomology, but generic deformations.

\begin{table}[ht]
\begin{center}
\begin{tabular}{rlllrrrr}
&Algebra&&Codifferential&$h_0$&$h_1$&$h_2$&$h_3$\\
\hline\\
&$d_{10}(p:q:r)$&=&$p\psa{1,4}1+\psa{2,4}1+q\psa{2,4}2+r\psa{3,4}3$&$0|0$&$2|0$&$0|2$&$0|0$\\
&$d_{10}(p:q:p-q)$&=&$p\psa{1,4}1+\psa{2,4}1+q\psa{2,4}2+(p-q)\psa{3,4}3$&$0|0$&$2|0$&$0|3$&$1|0$\\
&$d_{10}(p:0:q)$&=&$p\psa{1,4}1+\psa{2,4}1+q\psa{3,4}3$&$1|0$&$2|1$&$2|2$&$2|2$\\
&$d_{10}(p:2p:q)$&=&$p\psa{1,4}1+\psa{2,4}1+2p\psa{2,4}2+q\psa{3,4}3$&$0|0$&$2|0$&$1|2$&$0|1$\\
&$d_{10}(p:q:0)$&=&$p\psa{1,4}1+\psa{2,4}1+q\psa{2,4}2$&$0|1$&$3|0$&$0|3$&$1|0$\\
&$d_{10}(p:q:-2p)$&=&$p\psa{1,4}1+\psa{2,4}1+q\psa{2,4}2-2p\psa{3,4}3$&$0|0$&$2|0$&$0|3$&$1|0$\\
&$d_{10}(p:q:-p-q)$&=&$p\psa{1,4}1+\psa{2,4}1+q\psa{2,4}2-(p+q)\psa{3,4}3$&$0|0$&$2|0$&$0|3$&$1|0$\\
&$d_{10}(p:3p:q)$&=&$p\psa{1,4}1+\psa{2,4}1+3p\psa{2,4}2+q\psa{3,4}3$&$0|0$&$2|0$&$0|2$&$1|0$\\
&$d_{10}(p:q:2p-q)$&=&$p\psa{1,4}1+\psa{2,4}1+q\psa{2,4}2+(2p-q)\psa{3,4}3$&$0|0$&$2|0$&$0|2$&$0|1$\\
&$d_{10}(p:q:p)$&=&$p\psa{1,4}1+\psa{2,4}1+q\psa{2,4}2+p\psa{3,4}3$&$0|0$&$2|0$&$0|2$&$0|1$\\
&$d_{10}(p:-p:q)$&=&$p\psa{1,4}1+\psa{2,4}1-p\psa{2,4}2+q\psa{3,4}3$&$0|0$&$2|0$&$0|2$&$2|0$\\
&$d_{10}(p:q:-3p)$&=&$p\psa{1,4}1+\psa{2,4}1+q\psa{2,4}2-3p\psa{3,4}3$&$0|0$&$2|0$&$0|2$&$0|1$\\
&$d_{10}(p:q:-(2p+q))$&=&$p\psa{1,4}1+\psa{2,4}1+q\psa{2,4}2-(2p+q)\psa{3,4}3$&$0|0$&$2|0$&$0|2$&$0|1$\\
\hline\\
\end{tabular}
\end{center}
\caption{\protect Cohomology of subfamilies of $d_{10}(p:q:r)$}
\label{Table ZZ}
\end{table}

The special point $d_{10}(1:0:-1)$ has additional jump deformations to $d_3$, $d_5(0:1)$, $d_6(1:0)$ and deforms in a nbd of $d_5(0:1)$ and $d_6(1:0)$.
The special point $d_{10}(1:2:-1)$ has a jump deformation to $d_5(1:-1)$, $d_{10}(1:0:1)$ has a jump to $d_5(1:1)$, $d_{10}(1:2:1)$ has a jump  to

\begin{table}[ht]
\begin{center}
\begin{tabular}{rlllrrrr}
&Algebra&&Codifferential&$h_0$&$h_1$&$h_2$&$h_3$\\
\hline\\
&$d_{10}(p:q:r)$&=&$p\psa{1,4}1+\psa{2,4}1+q\psa{2,4}2+r\psa{3,4}3$&$0|0$&$2|0$&$0|2$&$0|0$\\
&$d_{10}(1:0:-1)$&=&$\psa{1,4}1+\psa{2,4}1-\psa{3,4}3$&$1|0$&$2|3$&$4|4$&$4|4$\\
&$d_{10}(1:2:-1)$&=&$\psa{1,4}1+\psa{2,4}1+2\psa{2,4}2-\psa{3,4}3$&$0|0$&$2|2$&$3|3$&$1|1$\\
&$d_{10}(1:0:1)$&=&$\psa{1,4}1+\psa{2,4}1+\psa{3,4}3$&$1|0$&$2|1$&$2|3$&$3|3$\\
&$d_{10}(1:2:1)$&=&$\psa{1,4}1+\psa{2,4}1+2\psa{2,4}2+\psa{3,4}3$&$0|0$&$2|0$&$1|3$&$1|2$\\
&$d_{10}(1:-1:2)$&=&$\psa{1,4}1+\psa{2,4}1-\psa{2,4}2+2\psa{3,4}3$&$0|0$&$2|0$&$0|4$&$4|0$\\
&$d_{10}(1:3:-2)$&=&$\psa{1,4}1+\psa{2,4}1+3\psa{2,4}2-2\psa{3,4}3$&$0|0$&$2|0$&$0|4$&$3|0$\\
&$d_{10}(1:0:0)$&=&$\psa{1,4}1+\psa{2,4}1$&$1|1$&$3|3$&$4|4$&$4|4$\\
&$d_{10}(1:0:-2)$&=&$\psa{1,4}1+\psa{2,4}1-2\psa{3,4}3$&$1|0$&$2|1$&$2|3$&$3|3$\\
&$d_{10}(1:2:0)$&=&$\psa{1,4}1+\psa{2,4}1+2\psa{2,4}2$&$0|1$&$3|0$&$1|3$&$1|2$\\
&$d_{10}(1:2:-2)$&=&$\psa{1,4}1+\psa{2,4}1+2\psa{2,4}2-2\psa{3,4}3$&$0|0$&$2|2$&$3|3$&$1|1$\\
&$d_{10}(1:2:-3)$&=&$\psa{1,4}1+\psa{2,4}1+2\psa{2,4}2-3\psa{3,4}3$&$0|0$&$2|0$&$1|3$&$1|2$\\
&$d_{10}(1:2:-4)$&=&$\psa{1,4}1+\psa{2,4}1+2\psa{2,4}2-4\psa{3,4}3$&$0|0$&$2|0$&$1|3$&$1|2$\\
&$d_{10}(1:-1:0)$&=&$\psa{1,4}1+\psa{2,4}1-\psa{2,4}2$&$0|1$&$3|0$&$0|4$&$4|0$\\
&$d_{10}(1:1:0)$&=&$\psa{1,4}1+\psa{2,4}1+\psa{2,4}2$&$0|1$&$3|0$&$0|3$&$1|0$\\
&$d_{10}(1:0:2)$&=&$\psa{1,4}1+\psa{2,4}1+2\psa{3,4}3$&$1|0$&$2|1$&$2|2$&$2|3$\\
&$d_{10}(1:0:-3)$&=&$\psa{1,4}1+\psa{2,4}1-3\psa{3,4}3$&$1|0$&$2|1$&$2|2$&$2|3$\\
&$d_{10}(1:3:-1)$&=&$\psa{1,4}1+\psa{2,4}1+3\psa{2,4}2-\psa{3,4}3$&$0|0$&$2|2$&$2|2$&$1|1$\\
&$d_{10}(1:3:5)$&=&$\psa{1,4}1+\psa{2,4}1+3\psa{2,4}2+5\psa{3,4}3$&$0|0$&$2|0$&$0|2$&$1|1$\\
&$d_{10}(1:3:1)$&=&$\psa{1,4}1+\psa{2,4}1+3\psa{2,4}2+\psa{3,4}3$&$0|0$&$2|0$&$0|2$&$1|1$\\
&$d_{10}(1:3:3)$&=&$\psa{1,4}1+\psa{2,4}1+3\psa{2,4}2+3\psa{3,4}3$&$0|0$&$2|0$&$0|2$&$1|1$\\
&$d_{10}(1:3:-3)$&=&$\psa{1,4}1+\psa{2,4}1+3\psa{2,4}2-3\psa{3,4}3$&$0|0$&$2|2$&$2|2$&$1|1$\\
&$d_{10}(1:3:-5)$&=&$\psa{1,4}1+\psa{2,4}1+3\psa{2,4}2-5\psa{3,4}3$&$0|0$&$2|0$&$0|2$&$1|1$\\
&$d_{10}(1:3:-7)$&=&$\psa{1,4}1+\psa{2,4}1+3\psa{2,4}2-7\psa{3,4}3$&$0|0$&$2|0$&$0|2$&$1|1$\\
&$d_{10}(1:3:-9)$&=&$\psa{1,4}1+\psa{2,4}1+3\psa{2,4}2-9\psa{3,4}3$&$0|0$&$2|0$&$0|2$&$1|1$\\
&$d_{10}(1:-1:3)$&=&$\psa{1,4}1+\psa{2,4}1-\psa{2,4}2+3\psa{3,4}3$&$0|0$&$2|0$&$0|2$&$2|2$\\
&$d_{10}(1:-3:5)$&=&$\psa{1,4}1+\psa{2,4}1-3\psa{2,4}2+5\psa{3,4}3$&$0|0$&$2|0$&$0|2$&$0|2$\\
&$d_{10}(1:5:-3)$&=&$\psa{1,4}1+\psa{2,4}1+5\psa{2,4}2-3\psa{3,4}3$&$0|0$&$2|0$&$0|2$&$0|2$\\
&$d_{10}(1:-1:1)$&=&$\psa{1,4}1+\psa{2,4}1-\psa{2,4}2+\psa{3,4}3$&$0|0$&$2|2$&$2|2$&$2|2$\\
&$d_{10}(1:-3:-3)$&=&$\psa{1,4}1+\psa{2,4}1-3\psa{2,4}2-3\psa{3,4}3$&$0|0$&$2|0$&$0|2$&$0|2$\\
&$d_{10}(1:-3:1)$&=&$\psa{1,4}1+\psa{2,4}1-3\psa{2,4}2+\psa{3,4}3$&$0|0$&$2|0$&$0|2$&$0|2$\\
&$d_{10}(0:0:0)$&=&$\psa{2,4}1$&$1|1$&$4|3$&$5|6$&$6|6$\\
\hline\\
\end{tabular}
\end{center}
\caption{\protect }
\label{Table W }
\end{table}

$d_5(2:1)$,
as well as deformations in a nbd of these points.  The special point $d_{10}(1:-1:2)$ has jumps to $d_5(1:2)$ and $d_7(1:-1)$, while $d_{10}(1:3:-2)$ has
jumps to $d_5(1:-2)$ and $d_7(1:3)$. The algebra $d_{10}(1:0:0)$ has jumps to $d_1$ and $d_7(0:1)$, and $d_{10}(1:0:-2)$ has a jump to $d_7(1:0)$, $d_{10}(1:2:-2)$ and $d_{10}(1:2:-4)$ have
jumps to $d_7(1:2)$.  The algebra $d_{10}(1:2:0)$ jumps to $d_1$, $d_{10}(1:2:-3)$ has a jump to $d_6(1:2)$, while $d_{10}(1:-1:0)$ has jumps to $d_1$
and $d_6(1:-1)$ and $d_{10}(1:1:0)$ has jumps to $d_1$ and $d_5(1:0)$.

The rest of the special points only have non generic cohomology, except $d_{10}(0:0:0)$ which jumps to $d_1$, $d_3$, $d_5(x:y)$ except $(0:0)$,
$d_6(x:y)$ except $(1:1)$ and $(0:0)$, $d_7(x:y)$ and $d_{10}(x:y:z)$ except $(0:0:0)$. Note this is exactly the deformation pattern for $d_5(0:0)$,
which is necessary as $d_5(0:0)$ and $d_{10}(0:0:0)$ are isomorphic algebras.

\section{The Moduli Space of Complex $1|3$-dimensional Lie Superalgebras}

\subsection{Construction of the moduli space of $1|3$-dimensional algebras}
Consider a $0|1$-dimensional vector space $W=\langle v_4\rangle$ and a $3|0$-dimensional vector space $M$.  There is no nontrivial $3|0$-dimensional codifferential, so $mu$ must vanish.  We have $\lambda=\sum_{i,j}\psa{j,4}ia_{i,j}$, which is given by a $3\times 3$ matrix
$A=\left[\begin{array}{ccc}a_{11}&a_{12}&a_{13}\\a_{21}&a_{22}&a_{23}\\a_{31}&a_{32}&a_{33}\end{array}\right]$. Moreover, $\psi$ also must vanish. Thus we
easily see that the algebras arising in this manner are given by the Jordan decomposition of the matrices $A$. This gives us three cases to consider,
given by the three matrices below:
\begin{equation*}
\left[\begin{array}{ccc}p&1&0\\0&q&1\\0&0&r\end{array}\right], \quad \left[\begin{array}{ccc}p&0&0\\0&p&1\\0&0&q\end{array}\right],\quad
\left[\begin{array}{ccc}1&0&0\\0&1&0\\0&0&1\end{array}\right].
\end{equation*}
The first matrix corresponds to the codifferential $d_1(p:q:r)$, the second to $d_2(p:q)$ and the third to $d_3$. There is a fourth case of Jordan
decomposition, given by the zero matrix, but that gives the trivial algebra, which we don't list in our table explicitly.

The other possible decomposition is given by the $1|0$-dimensional space $W=\langle v_3\rangle$ and
the $2|1$-dimensional space $M=\langle v_1,v_2,v_4\rangle$. There are 4 nontrivial possibilities for $\mu$, given by the elements in the $2|1$-dimensional
moduli space already described in this paper.  The $\lambda$ term is of the form $\lambda=\psa{1,3}4a_{31}+\psa{2,3}4a_{32}+\psa{3,4}1a_{1,3}+\psa{3,4}2a_{23}$,
while $\psi=\psa{3,3}4c$ and $\beta=\pha{3}1b_1+\pha{3}2b_2$.

The first case is $\mu=\psa{1,4}1p+\psa{2,4}1+\psa{2,4}2q$. The compatibility condition forces $a_{31}=a_{32}=0$.  If neither $p$ nor $q$ vanish, then
by adding an appropriate $[\mu,\beta]$ term, we could eliminate $\lambda$, but in general,  we can at least eliminate the $a_{13}$ term. Thus we reduce
to the case $\lambda=\psa{3,4}2a_{23}$. The MC equation forces $\psi=0$.  By applying an element of $G_{\delta,\mu}$ we can reduce to the case
where $a_{23}$ is either 1 or 0.  When $a_{23}=1$, we obtain the algebra $d_1(p: q: 0)$.  The case when $a_{23}=0$ is a bit more complex. When neither
$p$ nor $q$ vanishes, we obtain $d_1(p: q: 0)$, but when $q=0$ we obtain $d_2(0:p)$ and similarly when $p=0$ we obtain $d_2(0:q)$.

For $\mu=\psa{1,4}1+\psa{2,4}2$, the compatibility condition forces $a_{31}=a_{32}=0$, and then by adding a $[\mu,\beta]$ term, we can make $\lambda$ vanish.
The MC equation also forces $\psi$ to vanish.  Thus we get the codifferential $\mu$ which is $d_2(1: 0)$ in our list.

For $\mu=\psa{1,2}4+2\psa{1,1}4$, applying the compatibility condition and taking into account the addition of a coboundary term, $\lambda$ can be
made to vanish.  This time, the MC condition does not force $\psi=0$, but we can assume $c=1$, which gives $d_4$, or $c=0$, which gives $d_5$.

\begin{table}[ht]
\begin{center}
\begin{tabular}{rlllrrrr}
&Algebra&&Codifferential&$h_0$&$h_1$&$h_2$&$h_3$\\
\hline\\
&$d_{1}(p:q:r)$&=&$p\psa{1,4}1+\psa{2,4}1+q\psa{2,4}2+\psa{3,4}2+r\psa{3,4}3$&$0|0$&$2|0$&$0|2$&$0|0$\\
&$d_{2}(p:q)$&=&$p\psa{1,4}1+p\psa{2,4}2+\psa{3,4}2+q\psa{3,4}3$&$0|0$&$4|0$&$0|4$&$0|0$\\
&$d_{3}$&=&$\psa{1,4}1+\psa{2,4}2+\psa{3,4}3$&$0|0$&$8|0$&$0|8$&$0|0$\\
&$d_{4}$&=&$4\psa{1,1}4+\psa{2,3}4$&$0|1$&$4|0$&$8|0$&$12|0$\\
&$d_{5}$&=&$\psa{1,2}4+4\psa{1,1}4$&$1|1$&$5|1$&$9|1$&$13|1$\\
&$d_{6}$&=&$4\psa{1,1}4$&$2|1$&$7|2$&$12|3$&$17|4$\\
\hline
\end{tabular}
\end{center}
\caption{\protect Cohomology of $1|3$-Dimensional Complex Lie
Algebras}\label{Table 9}
\end{table}

For $\mu=2\psa{1,1}4$, we reduce to the case where only $a_{32}$ does not vanish, so $\lambda=\psa{3,4}2a_{32}$. We can assume that $a_{32}$ is either 1 or 0.
When $a_{32}=1$, the MC condition does not impose any restriction on $\psi$, nor does the cocycle condition. Nevertheless, independently of the value of $c$.
the algebra is isomorphic to $d_4$.  When $a_{32}=0$, we again don't get any restriction on $\psi$, but can take $c=1$ or $c=0$. The first gives $d_5$, while
the second gives $d_6$.

Finally, when $\mu=0$, the first restriction on $\lambda$ comes from the $MC$ equation, which forces either $a_{13}=a_{23}=0$, or $a_{31}=a_{32}=0$.
For the first case, we consider the action of $G_{\delta,\mu}$ on $\lambda$, which turns out to be equivalent to the action of $\GL(2,\C)$ on the
vector $(a_{31},a_{32})$. This action gives us only two equivalence classes, that of $(1,0)$ and $(0,0)$. The first class corresponds to
$a_{31}=1$ and $a_{32}=0$. There is no restriction on $\psi$, but independently of the value of $c$, we obtain $d_5$. Similarly, when both $a_{31}$ and
$a_{32}$ vanish, then we obtain $d_6$.   Finally, if $a_{31}=a_{32}=0$, then we get a similar action of $\GL(2,\C)$ on the vector $(a_{13},a_{23})$, again
with two equivalence classes. The nonvanishing class $a_{13}=1$ and $a_{23}=0$ forces $\psi=0$ by the cocycle condition, and we obtain $d_2(0:0)$.
The vanishing class gives $\lambda=0$, and there is no restriction on $\psi$.  When $c\ne0$ we obtain $d_6$. The case $c=0$ gives us  the trivial
algebra.  This completes the description of the moduli space.

\subsection{Deformations of the $1|3$-dimensional algebras}

\begin{table}[ht]
\begin{center}
\begin{tabular}{rlllrrrr}
&Algebra&&Codifferential&$h_0$&$h_1$&$h_2$&$h_3$\\
\hline\\
&$d_{1}(p:q:r)$&=&$p\psa{1,4}1+\psa{2,4}1+q\psa{2,4}2$\\
&&&$+\psa{3,4}2+r\psa{3,4}3$&$0|0$&$2|0$&$0|2$&$0|0$\\
&$d_{1}(p:q:0)$&=&$p\psa{1,4}1+\psa{2,4}1+q\psa{2,4}2+\psa{3,4}2$&$1|0$&$2|1$&$2|2$&$2|2$\\
&$d_{1}(p:2p:q)$&=&$p\psa{1,4}1+\psa{2,4}1+2p\psa{2,4}2$\\
&&&$+\psa{3,4}2+q\psa{3,4}3$&$0|0$&$2|0$&$1|2$&$0|1$\\
&$d_{1}(p:q:p+q)$&=&$p\psa{1,4}1+\psa{2,4}1+q\psa{2,4}2$\\
&&&$+\psa{3,4}2+(p+q)\psa{3,4}3$&$0|0$&$2|0$&$1|2$&$0|1$\\
&$d_{1}(p:3p:q)$&=&$p\psa{1,4}1+\psa{2,4}1+3p\psa{2,4}2$\\
&&&$+\psa{3,4}2+q\psa{3,4}3$&$0|0$&$2|0$&$0|2$&$1|0$\\
&$d_{1}(p:-p:q)$&=&$p\psa{1,4}1+\psa{2,4}1-p\psa{2,4}2$\\
&&&$+\psa{3,4}2+q\psa{3,4}3$&$0|0$&$2|0$&$0|2$&$2|0$\\
&$d_{1}(p:q:-2p+q)$&=&$p\psa{1,4}1+\psa{2,4}1+q\psa{2,4}2$\\
&&&$+\psa{3,4}2+(-2p+q)\psa{3,4}3$&$0|0$&$2|0$&$0|2$&$1|0$\\
&$d_{1}(p:q:2p+q)$&=&$p\psa{1,4}1+\psa{2,4}1+q\psa{2,4}2$\\
&&&$+\psa{3,4}2+(2p+q)\psa{3,4}3$&$0|0$&$2|0$&$0|2$&$1|0$\\
\hline\\
\end{tabular}
\end{center}
\caption{\protect Cohomology of special subfamilies of the family $d_1(p:q:r)$}
\label{Table U}
\end{table}

The family of algebras $d_1(p:q:r)$ is parametrized by $\P^2/\Sigma_3$, where $\Sigma_3$ acts by permuting the coordinates.
There are a lot of special subfamilies and

\begin{table}[ht]
\begin{center}
\begin{tabular}{rlllrrrr}
&Algebra&&Codifferential&$h_0$&$h_1$&$h_2$&$h_3$\\
\hline\\
&$d_{1}(1:-1:0)$&=&$\psa{1,4}1+\psa{2,4}1-\psa{2,4}2+\psa{3,4}2$&$1|0$&$2|1$&$3|2$&$4|3$\\
&$d_{1}(1:2:0)$&=&$\psa{1,4}1+\psa{2,4}1+2\psa{2,4}2+\psa{3,4}2$&$1|0$&$2|1$&$3|2$&$3|3$\\
&$d_{1}(1:2:4)$&=&$\psa{1,4}1+\psa{2,4}1+2\psa{2,4}2+\psa{3,4}2+4\psa{3,4}3$&$0|0$&$2|0$&$2|2$&$1|2$\\
&$d_{1}(1:2:3)$&=&$\psa{1,4}1+\psa{2,4}1+2\psa{2,4}2+\psa{3,4}2+3\psa{3,4}3$&$0|0$&$2|0$&$2|2$&$1|2$\\
&$d_{1}(1:-1:2)$&=&$\psa{1,4}1+\psa{2,4}1-\psa{2,4}2+\psa{3,4}2+2\psa{3,4}3$&$0|0$&$2|0$&$2|2$&$2|2$\\
&$d_{1}(1:-3:4)$&=&$\psa{1,4}1+\psa{2,4}1-3\psa{2,4}2+\psa{3,4}2+4\psa{3,4}3$&$0|0$&$2|0$&$1|2$&$0|1$\\
&$d_{1}(3:-4:-8)$&=&$3\psa{1,4}1+\psa{2,4}1-4\psa{2,4}2+\psa{3,4}2-8\psa{3,4}3$&$0|0$&$2|0$&$1|2$&$0|1$\\
&$d_{1}(3:-4:7)$&=&$3\psa{1,4}1+\psa{2,4}1-4\psa{2,4}2+\psa{3,4}2+7\psa{3,4}3$&$0|0$&$2|0$&$1|2$&$0|1$\\
&$d_{1}(1:-2:-4)$&=&$\psa{1,4}1+\psa{2,4}1-2\psa{2,4}2+\psa{3,4}2-4\psa{3,4}3$&$0|0$&$2|0$&$1|2$&$1|1$\\
&$d_{1}(1:0:0)$&=&$\psa{1,4}1+\psa{2,4}1+\psa{3,4}2$&$1|0$&$2|1$&$2|2$&$2|2$\\
&$d_{1}(1:3:0)$&=&$\psa{1,4}1+\psa{2,4}1+3\psa{2,4}2+\psa{3,4}2$&$1|0$&$2|1$&$2|2$&$3|2$\\
&$d_{1}(1:-1:3)$&=&$\psa{1,4}1+\psa{2,4}1-\psa{2,4}2+\psa{3,4}2+3\psa{3,4}3$&$0|0$&$2|0$&$0|2$&$4|0$\\
&$d_{1}(1:-3:3)$&=&$\psa{1,4}1+\psa{2,4}1-3\psa{2,4}2+\psa{3,4}2+3\psa{3,4}3$&$0|0$&$2|0$&$0|2$&$3|0$\\
&$d_{1}(1:3:7)$&=&$\psa{1,4}1+\psa{2,4}1+3\psa{2,4}2+\psa{3,4}2+7\psa{3,4}3$&$0|0$&$2|0$&$0|2$&$2|0$\\
&$d_{1}(1:3:9)$&=&$\psa{1,4}1+\psa{2,4}1+3\psa{2,4}2+\psa{3,4}2+9\psa{3,4}3$&$0|0$&$2|0$&$0|2$&$2|0$\\
&$d_{1}(1:3:5)$&=&$\psa{1,4}1+\psa{2,4}1+3\psa{2,4}2+\psa{3,4}2+5\psa{3,4}3$&$0|0$&$2|0$&$0|2$&$2|0$\\
&$d_{1}(1:3:-5)$&=&$\psa{1,4}1+\psa{2,4}1+3\psa{2,4}2+\psa{3,4}2-5\psa{3,4}3$&$0|0$&$2|0$&$0|2$&$2|0$\\
&$d_{1}(1:3:-7)$&=&$\psa{1,4}1+\psa{2,4}1+3\psa{2,4}2+\psa{3,4}2-7\psa{3,4}3$&$0|0$&$2|0$&$0|2$&$1|0$\\
&$d_{1}(1:-3:-5)$&=&$\psa{1,4}1+\psa{2,4}1-3\psa{2,4}2+\psa{3,4}2-5\psa{3,4}3$&$0|0$&$2|0$&$0|2$&$2|0$\\
&$d_{1}(1:-2:0)$&=&$\psa{1,4}1+\psa{2,4}1-2\psa{2,4}2+\psa{3,4}2$&$1|0$&$2|1$&$2|2$&$3|2$\\
&$d_{1}(1:-7:0)$&=&$\psa{1,4}1+\psa{2,4}1-7\psa{2,4}2+\psa{3,4}2$&$1|0$&$2|1$&$2|2$&$2|2$\\
&$d_{1}(3:-5:-7)$&=&$3\psa{1,4}1+\psa{2,4}1-5\psa{2,4}2+\psa{3,4}2-7\psa{3,4}3$&$0|0$&$2|0$&$0|2$&$1|0$\\
&$d_{1}(1:5:-9)$&=&$\psa{1,4}1+\psa{2,4}1+5\psa{2,4}2+\psa{3,4}2-9\psa{3,4}3$&$0|0$&$2|0$&$0|2$&$1|0$\\
&$d_{1}(1:3:-13)$&=&$\psa{1,4}1+\psa{2,4}1+3\psa{2,4}2+\psa{3,4}2-13\psa{3,4}3$&$0|0$&$2|0$&$0|2$&$1|0$\\
&$d_{1}(0:0:0)$&=&$\psa{2,4}1+\psa{3,4}2$&$1|0$&$3|1$&$4|3$&$5|4$\\
\hline\\
\end{tabular}
\end{center}
\caption{\protect Cohomology of special points in the family $d_1(p:q:r)$}
\label{Table T}
\end{table}

special points, for which the cohomology is not generic.
 However, none of these
special cases, except the generic point $d_1(0:0:0)$ give rise to any extra deformations, which makes sense, since $d_1(p:q:r)$ is the
first element in our list, so shouldn't have any extra deformations.  We do have a 2-parameter family of deformations, as
elements deform along the family.   The generic element $d_1(0:0:0)$ has jump deformations to element in the family except itself.

The family $d_2(p:q)$ is parametrized by $\P^1$, with no action of a symmetric group. Generically, an element $d_2(p:q)$ has
a jump deformation

\begin{table}[ht]
\begin{center}
\begin{tabular}{rlllrrrr}
&Algebra&&Codifferential&$h_0$&$h_1$&$h_2$&$h_3$\\
\hline\\
&$d_{2}(p:q)$&=&$p\psa{1,4}1+p\psa{2,4}2+\psa{3,4}2+q\psa{3,4}3$&$0|0$&$4|0$&$0|4$&$0|0$\\
&$d_{2}(0:1)$&=&$\psa{3,4}2+\psa{3,4}3$&$2|0$&$4|2$&$6|4$&$8|6$\\
&$d_{2}(1:0)$&=&$\psa{1,4}1+\psa{2,4}2+\psa{3,4}2$&$1|0$&$4|1$&$4|4$&$4|4$\\
&$d_{2}(1:2)$&=&$\psa{1,4}1+\psa{2,4}2+\psa{3,4}2+2\psa{3,4}3$&$0|0$&$4|0$&$3|4$&$0|3$\\
&$d_{2}(2:1)$&=&$2\psa{1,4}1+2\psa{2,4}2+\psa{3,4}2+\psa{3,4}3$&$0|0$&$4|0$&$2|4$&$0|2$\\
&$d_{2}(1:3)$&=&$\psa{1,4}1+\psa{2,4}2+\psa{3,4}2+3\psa{3,4}3$&$0|0$&$4|0$&$0|4$&$4|0$\\
&$d_{2}(3:1)$&=&$3\psa{1,4}1+3\psa{2,4}2+\psa{3,4}2+\psa{3,4}3$&$0|0$&$4|0$&$0|4$&$2|0$\\
&$d_{2}(1:-1)$&=&$\psa{1,4}1+\psa{2,4}2+\psa{3,4}2-\psa{3,4}3$&$0|0$&$4|0$&$0|4$&$6|0$\\
&$d_{2}(0:0)$&=&$\psa{3,4}2$&$2|0$&$5|2$&$8|5$&$11|8$\\
\hline\\
\end{tabular}
\end{center}
\caption{\protect Cohomology of the special points in the family $d_2(p:q)$}
\label{Table S}
\end{table}

to $d_1(p:p:q)$ and smooth deformations in nbds of $d_1(p:p:q)$ as well as $d_2(p:q)$. Again, there are
special points, but no extra deformations, because the elements already deform to everything they could. The exception is
$d_2(0:0)$, which has jump deformations to $d_1(x: y:z)$ for all $(x:y:z)$ and $d_2(x:y)$ for all $(x:y)$ except $(0:0)$.

The algebra $d_3$ has jump deformations to $d_1(1:1:1)$ and $d_2(1:1)$ as well as smooth deformations in nbds of these points.
The description of the deformation picture of the first 3 algebras corresponds exactly to the description of the moduli space
of $3\times 3$ matrices with the action given by conjugation by an element in $\GL(3,\C)$ and multiplication by a nonzero complex
number.  This picture arises in many of the moduli spaces of different algebraic objects.

The algebra $d_4$ is rigid, even though the dimension of $H^2$ is 8, because it is really $8|0$-dimensional, and only odd elements
of $H^2$ contribute to deformations.  The algebra $d_5$ has a jump deformation to $d_4$, and $d_6$  has jump deformations to $d_4$
and $d_5$.  This completes the description of the deformations of the elements in the moduli space of $1|3$-dimensional complex
Lie superalgebras.

\bibliographystyle{amsplain}

\providecommand{\bysame}{\leavevmode\hbox to3em{\hrulefill}\thinspace}
\providecommand{\MR}{\relax\ifhmode\unskip\space\fi MR }
\providecommand{\MRhref}[2]{%
  \href{http://www.ams.org/mathscinet-getitem?mr=#1}{#2}
}
\providecommand{\href}[2]{#2}

\end{document}